\documentclass[11pt]{amsart}
\usepackage{fullpage}
\usepackage[utf8]{inputenc}
\usepackage[T1]{fontenc}
\usepackage{graphicx}
\usepackage{amsmath}
\usepackage{amsfonts}
\usepackage{graphicx}

\newtheorem{theo}{Theorem}[section]
\newtheorem{cor}[theo]{Corollary}
\newtheorem{rem}[theo]{Remark}

\newtheorem{propo}[theo]{Proposition}
\newtheorem{lemme}[theo]{Lemma}

\newtheorem{ex}[theo]{Example}

\newtheorem{hyp}[theo]{Assumptions}
\newcommand{\E}{\mathbb{E}}
\newcommand{\R}{\mathbb{R}}
\newcommand{\PP}{\mathbb{P}}

\newcommand{\F}{\mathcal{F}}

\newcommand{\N}{\mathbb{N}}

\title[Invariant measure approximation for SPDE's]{Approximation of the invariant measure with an Euler scheme for Stochastic PDE's driven by Space-Time White Noise}
\author[C-E BREHIER]{Charles-Edouard BREHIER}
\address{ENS Cachan Bretagne - IRMAR, Universit\'e Rennes 1\\Avenue Robert Schumann\\F-35170 Bruz\\France}
\email{charles-edouard.brehier@bretagne.ens-cachan.fr}
\date{}

\begin{document}
\maketitle
\begin{abstract}

In this article, we consider a stochastic PDE of parabolic type, driven by a space-time white-noise, and its numerical discretization in time with a semi-implicit Euler scheme. When the nonlinearity is assumed to be bounded, then a dissipativity assumption is satisfied, which ensures that the SDPE admits a unique invariant probability measure, which is ergodic and strongly mixing - with exponential convergence to equilibrium. Considering test functions of class $\mathcal{C}^2$, bounded and with bounded derivatives, we prove that we can approximate this invariant measure using the numerical scheme, with order $1/2$ with respect to the time step.







\end{abstract}

\section{Introduction}

In this paper, we are interested in the discretization in time of the following stochastic reaction-diffusion equation
\begin{equation}\label{truePDE}
\frac{\partial y(t,\xi)}{\partial t}=\frac{\partial^2 y(t,\xi)}{\partial \xi^2}+g(\xi,y(t,\xi))+\frac{\partial \omega(t,\xi)}{\partial t},
\end{equation}
for $t\geq 0, \xi\in(0,1)$, with the initial condition $y(0,\xi)=y(\xi)$, and homogeneous Dirichlet boundary conditions $y(t,0)=y(t,1)=0$.
The stochastic perturbation $\frac{\partial \omega(t,\xi)}{\partial t}$ is a space-time white noise: 
the rigorous interpretation of \eqref{truePDE} is given by an abstract evolution equation \eqref{eqY} - in the sense of \cite{DaP-Z1} - in the Hilbert space $H=L^2(0,1)$, driven by a Wiener cylindrical process - see Section \ref{SectWiener}.

The numerical approximation of stochastic equations has been extensively studied during the last thirty years. First we recall that one can look at strong approximation results - when the trajectories of the continuous and the discrete-time processes are compared - or at weak approximation results - when the laws at a fixed time are compared. The simplest method in the case of SDEs is the Euler-Maruyama scheme; it is built as a straightforward extension of the well-known explicit Euler method for ODEs, which is of order $1$: if we consider in $\R^d$ a SDE with regular coefficients
$$dX_t=f(X_t)dt+\sigma(X_t)dB_t,X_0=x,$$
its numerical approximation is defined for a given step-size $\Delta t$ by
\begin{gather*}
X_0=x,\\
X_{n+1}=X_n+f(X_n)\Delta t+\sigma(X_n)(B_{t_{n+1}}-B_{t_{n}}).
\end{gather*}
Due to the regularity properties of the Brownian Motion, this scheme is in general only of order $1/2$ in the strong sense - i.e. for $n\Delta t\leq T$ we have $\E|X_{n\Delta t}-X_n|^2\leq C(T)\Delta t$ - while it is of order $1$ in the weak sense - i.e. for test functions $\phi$ of class $\mathcal{C}^3$, with bounded derivatives, we have a bound $|\E\phi(X_{n\Delta t})-\E\phi(X_n)|\leq C(T,\phi)\Delta t$. The idea for proving that weak order is $1$ and not $1/2$ is to consider the Kolmogorov equation associated with the process $X_t$, which is satified by the function $(t,x)\mapsto \E\phi(X(t,x))$ - see \cite{TalGen}, \cite{TalTub}. The books \cite{KloPla} and \cite{MilTre} - see also \cite{Mil_book} - contain various numerical schemes - like the well-known Milstein scheme, some implicit schemes, and methods based on stochastic Taylor expansions - with their order of convergence in both strong and weak senses.




Numerical methods for SPDEs like \eqref{truePDE} need discretization both in time and in space. For example, time discretization leads to explicit or implicit methods, while discretization in space can be done with finite difference or finite element methods. Basically, the result for space-time white noise driven equations is convergence with strong order $1/2$ in space and only $1/4$ in time - under some Courant-Friedrichs-Lewy conditions when necessary: see \cite{davie-gaines}, \cite{gyongy1}, \cite{gyongy2}, \cite{gyongy-nualart}, \cite{haus}, \cite{walsh}.

If we look at the abstract formulation of \eqref{truePDE} in the Hilbert space $H$, results have also been proved for the time discretization using semi-implicit Euler schemes: the strong order of convergence is $1/4$ - see \cite{prin2} - while the weak order of convergence is $1/2$ - see \cite{deb}. Moreover in \cite{deb-prin}, the authors have studied weak convergence in the case of linear equations when using a finite element method in space. We follow here the framework of \cite{deb}: we consider the stochastic evolution equation in $H$ and we use a time discretization with a semi-implicit Euler scheme, with no discretization in space.






In this work, we are interested in the behaviour of the weak convergence estimates when the final time $T$ goes to infinity: can we replace constants $C(T,\phi)$ by a constant $C(\phi)$ independent from $T$? Passing to the limit, we thus ask the more general following question: can we use a numerical scheme to approximate the invariant probability measure of the continuous time process - which is assumed to be unique, ergodic and with exponential convergence to equilibrium? The SDE case has been studied with various methods: in \cite{TalInv}, the weak error analysis is made by showing that the time derivatives of the solution of the Kolmogorov equation are exponentially decreasing in time; in \cite{HighMattStu}, some general conditions are given for the ergodicity of the numerical scheme, thanks to the theory of geometric ergodicity of Markov Chains. In \cite{MattStuTre}, the approximation result is shown thanks to the use of a Poisson equation.

In the case of SPDEs, general ergodicity results have only been obtained recently - see Section \ref{sectasymptu} - and the problem of approximation of invariant measures by numerical schemes has not been studied yet. 



We prove the following result:
\begin{theo}\label{th}
For any $0<\kappa<1/2$, $\tau_0>0$ and for any $\mathcal{C}_{b}^{2}$ function $\phi$, there exists a constant $C>0$ such that for any $m\geq 2$, $y\in H$ and $0<\tau\leq \tau_0$
$$|\E[\phi(Y(m\tau,y))]-\E[\phi(Y_{m}(\tau,y))]|\leq C(1+|y|^3)(((m-1)\tau)^{-1/2+\kappa}+1)\tau^{1/2-\kappa}.$$
\end{theo}
The continuous process $Y$ is defined by Equation \eqref{eqY} below, and the numerical approximation $(Y_m(\tau,y))$ with time step $\tau$ and initial condition $y$ is defined by \eqref{defYk}.

We notice that the right-hand side of the previous estimate contains a singularity when $m=1$, which is due to a lack of regularity of the infinite-dimensional processes - details are given in Sub-Section \ref{lackreg} below. However, if we look at the error at a fixed time $T=m\tau$, if $\tau$ is small enough we just need to change the constant $C=C(T)$; moreover since we are interested in the behaviour when $m$ goes to infinity, this term plays no role.


With the assumptions precised below, we show in Section \ref{sectasymptu} that the SPDE admits a unique invariant probability measure $\overline{\mu}$, which is ergodic and strongly mixing, with exponential convergence to equilibrium; nevertheless we can in general only show the existence, not the uniqueness, of invariant measures for the numerical approximation, and we can prove the following result:
\begin{cor}\label{corMesInv}
For any $0<\kappa<1/2$, $\tau_0>0$ and for any $\mathcal{C}_{b}^{2}$ function $\phi$, there exists constants $c>0$, $C>0$ such that for any $0<\tau\leq \tau_0$, any initial condition $y\in H$ and any $m\geq 1$
$$|\E[\phi(Y_{m}(\tau,y))]-\int_{H}\phi d\overline{\mu}|\leq C(1+|y|^3)(\frac{1}{m^{1/2-\kappa}}+\tau^{1/2-\kappa})+C(1+|y|^2)e^{-cm\tau}.$$
Moreover, if $\mu^\tau$ is an ergodic invariant probability measure of the numerical scheme $(Y_m(\tau,.))_{m\in \N}$, we have
$$|\int_{H}\phi d\mu^\tau-\int_{H}\phi d\overline{\mu}|\leq C\tau^{1/2-\kappa};$$
all the ergodic invariant measures of the numerical approximation are then close to the unique invariant probability measure of the continuous process.
\end{cor}


Up to our knowledge, this is the first result of this kind for SPDEs.




The key point, like in \cite{TalInv}, for obtaining bounds independent from the time $T=m\tau$, is to prove that derivatives of the solution $u$ of the underlying Kolomogorov equation - mentioned above - decrease exponentially in time: this is done in Section \ref{expdec}, using a coupling method.

Moreover, an essential tool in \cite{deb} is the use of Malliavin calculus and of an integration by parts formula in order to transform a stochastic integral - which is not regular enough in space in the infinite dimensional setting - into an expression which can be controlled; in Lemma \ref{lem5}, we see that the involved Malliavin derivatives may increase exponentially fast with respect to time. The solution is to separate the lack of regularity problem and this badly controlled growth with a decomposition of the interval into two parts: in the first one we need the integration by parts formula, and we can control the Malliavin derivatives, while in the other one we can directly give an appropriate bound. We also provide an improvement with respect to \cite{deb}: we can consider a more general nonlinear coefficient $G$ - like a Nemytskii operator, see Example \ref{exG}.

We also notice that some more general equations, with additive noise which is white in time but colored in space, can be studied with our method, with suitable assumptions on the coefficients. For a very smooth noise, the numerical analysis on finite time is easier to treat, but then ergodic properties and long-time behaviour require different techniques; this will be treated elsewhere.

The case of some equations with multiplicative noise which satisfy the Strong Feller Property is also covered by our technique of proof, but with some additional difficulties - see for instance the restrictive condition on the diffusion coefficient in \cite{deb}. Since all the necessary ideas are contained here and in \cite{deb}, we only focus on the additive noise case.

The paper is organized as follows: in Section \ref{sectass}, we precise the assumptions made on the coefficients of the equations, and we define the numerical method in Section \ref{sectdefnum}. In Section \ref{sectprel}, we give some bounds on the solutions of the continuous and discrete equations, and we study their asymptotic behaviour: existence of invariant measures, and uniqueness for the continuous equation, following from Proposition \ref{propoexpy1y2}. In Section \ref{sectproof}, we explain the proof of Theorem \ref{th}, and we give a proof of Corollary \ref{corMesInv}; more precisely, in Section \ref{expdec} we prove the exponential decreasing in time of the derivatives of the function $u$. Eventually Section \ref{sectestims} contains the proofs of the remaining estimates


\section{Notations and assumptions}\label{sectass}

Let $H$ be a separable Hilbert space, with norm denoted by $|.|_{H}$ or simply $|.|$. We consider equations of the form
\begin{equation}\label{eqY}
\begin{gathered}
dY(t,y)=(BY(t,y)+G(Y(t,y)))dt+dW(t)\\
Y(0,y)=y.
\end{gathered}
\end{equation}
In the next paragraphs, we explain the assumptions on the linear operator $B$ and on the nonlinear coefficient $G$; we also recall how the cylindrical Wiener process $W$ is defined, and how we can construct solutions of this equation.

\subsection{Test functions}
To quantify the weak approximation, we use test functions $\phi$ in the space $\mathcal{C}_{b}^{2}(H,\R)$ of functions from $H$ to $\R$ that are twice continuously differentiable, bounded, with first and second order bounded derivatives.

\begin{rem}\label{remident}
In the sequel, we often identify the first derivative $D\phi(x)\in \mathcal{L}(H,\R)$ with the gradient in the Hilbert space $H$, and the second derivative $D^2\phi(x)$ with a linear operator on $H$, via the formulas:
\begin{gather*}
<D\phi(x),h>=D\phi(x).h \text{ for every }h\in H\\
<D^2\phi(x).h,k>=D^2\phi(x).(h,k) \text{ for every }h,k\in H.
\end{gather*}
\end{rem}

In the sequel, we use the following notations:
\begin{gather*}
\|\Phi\|_{\infty}=\sup_{x\in H}|\Phi(x)|_{H}\\
\|\Phi\|_{1}=\sup_{x\in H}|D\Phi(x)|_{H}\\
\|\Phi\|_{2}=\sup_{x\in H}|D^2\Phi(x)|_{\mathcal{L}(H)}.
\end{gather*}

\subsection{Assumptions on the coefficients}

\subsubsection{The linear operator}

We denote by $\N=\left\{0,1,2,\ldots\right\}$ the set of nonnegative integers.

We suppose that the following properties are satisfied:
\begin{hyp}\label{hypB}
\begin{enumerate}
\item We assume that there exists a complete orthonormal system of elements of $H$ denoted by $(f_k)_{k\in \N}$, and a non-decreasing sequence of real positive numbers $(\mu_k)_{k\in \N}$ such that:
$$Bf_{k}=-\mu_{k}f_{k}\text{ for all } k\in\N.$$
\item The sequence $(\mu_k)$ goes to $+\infty$ and
$$
\sum_{k=0}^{+\infty}\frac{1}{\mu_{k}^{\alpha}}<+\infty\Leftrightarrow \alpha>1/2.$$
\end{enumerate}
\end{hyp}

The smallest eigenvalue of $-B$ is then $\mu_0$.

\begin{ex}\label{exampleB}
The equation \eqref{truePDE} enters in this framework: we can choose $B=\frac{d^2}{dx^2}$, with the domain $H^2(0,1)\cap H_{0}^{1}(0,1)\subset L^2(0,1)$ - corresponding to homogeneous Dirichlet boundary conditions. In this case for any $k\in\N$ $\mu_k=\pi^2 (k+1)^2$, and $f_k(\xi)=\sqrt{2}\sin((k+1)\pi\xi)$ - see \cite{Brezis}.
\end{ex}

For a $N\in\left\{1,2,\ldots\right\}$, we denote by $H_N$ the subspace of $H$ spanned by $f_0,\ldots,f_{N-1}$, and by $P_N$ the orthogonal projection of $H$ onto $H_N$.

The domain $D(B)$ of $B$ is equal to $D(B)=\left\{y=\sum_{k=0}^{+\infty}y_kf_{k}\in H, \sum_{k=0}^{+\infty}(\mu_k)^{2}|y_k|^2<+\infty\right\}$.
We can more generally define fractional powers of $-B$, for $b\in[0,1]$:
$$(-B)^b y=\sum_{k=0}^{\infty}\mu_{k}^{b}y_kf_k\in H,$$
with the domains
$$D(-B)^b=\left\{y=\sum_{k=0}^{+\infty}y_kf_{k}\in H,\quad |y|_{b}^{2}:=\sum_{k=0}^{+\infty}(\mu_k)^{2b}|y_k|^2<+\infty\right\}.$$

When $b\in[0,1]$, we can also define the spaces $D(-B)^{-b}$ and operators $(-B)^{-b}$, with norm denoted by $|.|_{-b}$; when $y=\sum_{k=0}^{+\infty}y_kf_{k}\in H$, we have $(-B)^{-b}y=\sum_{k=0}^{+\infty}\mu_{k}^{-b}y_kf_k$ and $|y|_{-b}^{2}:=\sum_{k=0}^{+\infty}(\mu_k)^{-2b}|y_k|^2$.

The semi-group $(e^{tB})_{t\geq 0}$ can be defined by the Hille-Yosida Theorem - see \cite{Brezis}. We use the following spectral formula: if $y=\sum_{k=0}^{+\infty}y_kf_k\in H$, then for any $t\geq 0$
$$e^{tB}y=\sum_{k=0}^{+\infty}e^{-\mu_k t}y_kf_k.$$

For any $t\geq 0$, $e^{tB}$ is a continuous linear operator in $H$, with operator norm $e^{-\mu_0 t}$. The semi-group $(e^{tB})$ is used to define the solution $Z(t)=e^{tB}z$ of the linear Cauchy problem
$$\frac{dZ(t)}{dt}=BZ(t)\quad \text{with} \quad Z(0)=z.$$

To define solutions of more general PDEs of parabolic type, we use mild formulation, and Duhamel principle.

This semi-group enjoys some smoothing properties that we often use in this work. Basically we need the following properties, which are easily proved using the above spectral properties.
\begin{propo}\label{proporegul}
Under Assumption \ref{hypB}, for any $\sigma\in[0,1]$, there exists $C_\sigma>0$ such that we have:
\begin{enumerate}
\item for any $t>0$ and $y\in H$,
$$|e^{tB}y|_{\sigma}\leq C_\sigma t^{-\sigma}e^{-\frac{\mu_0}{2}t}|y|_{H}.$$
\item for any $0<s<t$ and $y\in H$,
$$|e^{tB}y-e^{sB}y|_{H}\leq C_\sigma\frac{(t-s)^\sigma}{s^\sigma}e^{-\frac{\mu_0}{2}s}|y|_H.$$
\item for any $0<s<t$ and $y\in D(-A)^\sigma$,
$$|e^{tB}y-e^{sB}y|_{H}\leq C_\sigma(t-s)^\sigma e^{-\frac{\mu_0}{2}s}|y|_{\sigma}.$$
\end{enumerate}
\end{propo}

\subsubsection{The nonlinear operator}

The nonlinear operator $G$ is assumed to satisfy some general assumptions. In Example \ref{exG}, we give the two main kind of operators that can be used in our framework.

\begin{hyp}\label{hypG}
The function $G:H\rightarrow H$ is assumed to be bounded and Lipschitz continuous. We denote by $L_G$ the Lipschitz constant of $G$

We also define for each $N\geq 1$ a function $G_N:H_N\rightarrow H_N$, with $G_N(y)=P_NG(y)$ for any $y\in H_N$. We assume that each $G_N$ is twice differentiable, and that we have the following bounds on the derivatives, uniformly with respect to $N$:
\begin{itemize}
\item There exists a constant $C_1$ such that for any $N\geq 1$, $y\in H_N$ and $h\in H_N$
$$|DG_N(y).h|_{H}\leq C_1|h|_{H}.$$
\item There exists $\eta\in[0,1)$ and a constant $C_2$ such that for any $N\geq 1$, $y\in H_N$ and any $h,k\in H_N$ we have
$$|(-B)^\eta D^2G_N(y).(h,k)|\leq C_2|h|_H|k|_H.$$
\item Moreover, there exists a constant $C_3$ such that for any $N\geq 1$, $y\in H_N$ and any $h,k\in H_N$
$$|D^2G_N(y).(h,k)|\leq C_3|h|_{(-B)^\eta}|k|_{H}.$$
\end{itemize}
\end{hyp}

Since $G$ is bounded, the following property is easily satisfied:
\begin{propo}[Dissipativity]\label{propodiss}
There exist $c>0$ and $C>0$ such that for any $y\in D(B)$
\begin{equation}\label{hypdiss}
<By+G(y),y>\leq -c|y|^2+C.
\end{equation}
\end{propo}

We remark that we have uniform control with respect to the dimension $N$ of the bounds on $G_N$ and on its derivatives, and that \eqref{hypdiss} is also satisfied for $G_N$, with constants $c$ and $C$ independent from $N$.


\begin{ex}\label{exG}
We give some fundamental examples of nonlinearities for which the previous assumptions are satisfied:
\begin{itemize}
\item A function $G:H\rightarrow H$ of class $\mathcal{C}^2$, bounded and with bounded derivatives, 
fits in the framework, with the choice $\eta=0$.
\item The function $G$ can be a \textbf{Nemytskii} operator: let $g:(0,1)\times \R\rightarrow \R$ be a measurable, bounded, function such that for almost every $\xi\in(0,1)$ $g(\xi,.)$ is twice continuously differentiable, with uniformly bounded derivatives.
Then $G(y)$ is defined for every $y\in H=L^2(0,1)$ by
$$G(y)(\xi)=g(\xi,y(\xi)).$$
In general, such functions are not Fr\'echet differentiable, but only G\^ateaux differentiable, with the following expressions:
\begin{gather*}
[DG(y).h](\xi)=\frac{\partial g}{\partial y}(\xi,y(\xi))h(\xi)\\
[D^2G(y).(h,k)](\xi)=\frac{\partial^2 g}{\partial y^2}(\xi,y(\xi))h(\xi)k(\xi).
\end{gather*}
If $h$ and $k$ are only $L^2$ functions, $D^2G(y).(h,k)$ may only be $L^1$; however if $h$ or $k$ is $L^\infty$, it is $L^2$.
The conditions in Assumption \ref{hypG} are then satisfied as soon as there exists $\eta<1$ such that $D(-B)^\eta$ is continuously embedded into $L^\infty(0,1)$ - it is the case for $B$ given in Example \ref{exampleB}, with $\eta>1/4$. Then the finite dimensional spaces $H_N$ are subspaces of $L^\infty$, and differentiability can be shown.
\end{itemize}
\end{ex}

\subsection{The cylindrical Wiener process and stochastic integration in $H$}\label{SectWiener}

In this section, we recall the definition of the cylindrical Wiener process and of stochastic integral on a separable Hilbert space $H$ with norm $|.|_H$. For more details, see \cite{DaP-Z1}.

We first fix a filtered probability space $(\Omega,\mathcal{F},(\mathcal{F}_t)_{t\geq 0},\PP)$. A cylindrical Wiener process on $H$ is defined with two elements:
\begin{itemize}
\item a complete orthonormal system of $H$, denoted by  $(q_i)_{i\in I}$, where $I$ is a subset of $\N$;
\item a family $(\beta_i)_{i\in I}$ of independent real Wiener processes with respect to the filtration $((\mathcal{F}_t)_{t\geq 0})$;
\end{itemize}

then $W$ is defined by
\begin{equation}\label{defWiener}
W(t)=\sum_{i\in I}\beta_i(t)q_i.
\end{equation}

When $I$ is a finite set, we recover the usual definition of Wiener processes in the finite dimensional space $\R^{|I|}$. However the subject here is the study of some Stochastic Partial Differential Equations, so that in the sequel the underlying Hilbert space $H$ is infinite dimensional; for instance when $H=L^{2}(0,1)$, an example of complete orthonormal system is $(q_k)=(\sqrt{2}\sin(k\pi.))_{k\geq 1}$ - see Example \ref{exampleB}.

A fundamental remark is that the series in \eqref{defWiener} does not converge in $H$; but if a linear operator $\Psi:H\rightarrow K$ is Hilbert-Schmidt, then $\Psi W(t)$ converges in $L^2(\Omega,H)$ for any $t\geq 0$.

We recall that a bounded linear operator $\Psi:H\rightarrow K$ is said to be Hilbert-Schmidt when
$$|\Psi|_{\mathcal{L}_{2}(H,K)}^{2}:=\sum_{k=0}^{+\infty}|\Psi(q_k)|_{K}^{2}<+\infty,$$
where the definition is independent of the choice of the orthonormal basis $(q_k)$ of $H$.
The space of Hilbert-Schmidt operators from $H$ to $K$ is denoted $\mathcal{L}_{2}(H,K)$; endowed with the norm $|.|_{\mathcal{L}_{2}(H,K)}$ it is an Hilbert space.

The stochastic integral $\int_{0}^{t}\Psi(s)dW(s)$ is defined in $K$ for predictible processes $\Psi$ with values in $\mathcal{L}_2(H,K)$ such that $\int_{0}^{t}|\Psi(s)|_{\mathcal{L}_2(H,K)}^{2}ds<+\infty$ a.s; moreover when $\Psi\in L^2(\Omega\times[0,t];\mathcal{L}_2(H,K))$, the following two properties hold:
\begin{gather*}
\E|\int_{0}^{t}\Psi(s)dW(s)|_{K}^{2}=\E\int_{0}^{t}|\Psi(s)|_{\mathcal{L}_2(H,K)}^{2}ds \text{ (It\^o isometry),}\\
\E\int_{0}^{t}\Psi(s)dW(s)=0.
\end{gather*}
A generalization of It\^o formula also holds - see \cite{DaP-Z1}.

For instance, if $v=\sum_{k\in\N}v_kq_k\in H$, we can define $$<W(t),v>=\int_{0}^{t}<v,dW(s)>=\sum_{k\in\N}\beta_k(t)v_k;$$
we then have the following space-time white noise property
$$\E<W(t),v_1><W(s),v_2>=t\wedge s<v_1,v_2>.$$

Therefore to be able to integrate a process with respect to $W$ requires some strong properties on the integrand; in our SPDE setting, the Hilbert-Schmidt properties follow from the assumptions made on the linear coefficients of the equations.

Thanks to Assumption \ref{hypB}, it is easy to show that the following stochastic integral is well-defined in $H$, for any $t\geq 0$:
\begin{equation}\label{stoconvWB}
W^B(t)=\int_{0}^{t}e^{(t-s)B}dW(s).
\end{equation}

It is called a stochastic convolution, and it is the unique mild solution of
$$dZ(t)=BZ(t)dt+dW(t)\quad \text{with} \quad Z(0)=0.$$

Under the second condition of Assumption \ref{hypB}, there exists $\delta>0$ such that for any $t>0$ we have $\int_{0}^{t}\frac{1}{s^\delta}|e^{sB}|_{\mathcal{L}_{2}(H)}^{2}ds<+\infty$; it can then be proved that $W^B$ has continuous trajectories - via the \textit{factorization method}, see \cite{DaP-Z1} - and that for any $1\leq p<+\infty$
\begin{equation}\label{momWBp}
\E\sup_{t\geq 0}|W^B(t)|_{H}^{p}<+\infty.
\end{equation}

We can now define solutions to Equation \eqref{eqY}, thanks to the assumptions made on the coefficients: the following result is classical - see \cite{DaP-Z1}:
\begin{propo}
For every $T>0$, $y\in H$, the equation \eqref{eqY} admits a unique mild solution $Y\in \text{L}^2(\Omega,\mathcal{C}([0,T],H))$:
\begin{equation}\label{eqmild}
Y(t)=e^{tB}y+\int_{0}^{t}e^{(t-s)B}G(Y(s))ds+\int_{0}^{t}e^{(t-s)B}dW(s).
\end{equation}
\end{propo}

\section{Definition of the numerical scheme}\label{sectdefnum}

We now define the numerical approximation of $Y$: denoting by $\tau$ the time step, we have
\begin{gather*}
Y_{k+1}(\tau,y)=Y_{k}(\tau,y)+\tau BY_{k+1}(\tau,y)+\tau G(Y_k(\tau,y))+\sqrt{\tau}\chi_{k+1}\\
Y_{0}(\tau,y)=y,
\end{gather*}
where $\chi_{k+1}=\frac{1}{\sqrt{\tau}}(W((k+1)\tau)-W(k\tau))$.

To simplify the equations, we omit the dependence of $Y_k$ on the time-step $\tau$ and on the initial condition $y$.

This expression does not make sense in $H$. Defining $R_\tau=(I-\tau B)^{-1}$, this last equation can be replaced by
\begin{equation}\label{defYk}
Y_{k+1}=R_\tau Y_{k}+\tau R_\tau G(Y_k)+\sqrt{\tau}R_\tau\chi_{k+1},
\end{equation}
which is valid, since $R_\tau$ is a Hilbert-Schmidt operator on $H$.

\begin{rem}
Later, we often use the following expression for $Y_k$:
\begin{equation}\label{exprYk}
Y_k=R_{\tau}^{k}y+\tau\sum_{l=0}^{k-1}R_{\tau}^{k-l}G(Y_l)+\sqrt{\tau}\sum_{l=0}^{k-1}R_{\tau}^{k-l}\chi_{l+1}.
\end{equation}
The following expression is also useful:
\begin{equation}\label{exprsto}
\sqrt{\tau}\sum_{l=0}^{k-1}R_{\tau}^{k-l}\chi_{l+1}=\int_{0}^{t_k}R_{\tau}^{k-l_s}dW(s),
\end{equation}
where $l_s=\lfloor \frac{s}{\tau} \rfloor$ - with the notation $\lfloor . \rfloor$ for the integer part.
\end{rem}

We need the following technical estimate:
\begin{lemme}\label{lem6}
For any $0\leq \kappa\leq 1$ and $j\geq 1$,
$$|(-B)^{1-\kappa}R_{\tau}^{j}|_{\mathcal{L}(H)}\leq c\frac{1}{(j\tau)^{1-\kappa}}\frac{1}{(1+\mu_0\tau)^{j\kappa}}.$$
\end{lemme}

\underline{Proof}
For any $z\in H$,
\begin{align*}
|(-B)^{1-\kappa}R_{\tau}^{j}z|_{H}^{2}&=\sum_{i=0}^{+\infty}\mu_{i}^{2(1-\kappa)}\frac{1}{(1+\mu_i\tau)^{2j}}|z_i|^2\\
&=\frac{1}{(j\tau)^{2(1-\kappa)}}\sum_{i=0}^{+\infty}|z_i|^2\mu_{i}^{2(1-\kappa)}(j\tau)^{2(1-\kappa)}\frac{1}{(1+\mu_i\tau)^{2j(1-\kappa)}}\frac{1}{(1+\mu_i\tau)^{2j\kappa}}\\
&\leq \frac{1}{(j\tau)^{2(1-\kappa)}}\sum_{i=0}^{+\infty}\left(\frac{\mu_i j\tau}{1+\mu_i j\tau}\right)^{2(1-\kappa)}\frac{1}{(1+\mu_0\tau)^{2j\kappa}}|z_i|^2\\
&\leq c|z|_{H}^{2}\frac{1}{(j\tau)^{2(1-\kappa)}}\frac{1}{(1+\mu_0\tau)^{2j\kappa}}.\qed
\end{align*}


\section{Preliminary results}\label{sectprel}


We warn the reader that constants may vary from line to line during the proofs, and that in order to use lighter notations we usually forget to mention dependence on the parameters. We use the generic notation $C$ for such constants.

We fix the time step $\tau$, as well as $m\in \N$; we then introduce the notation $T=m\tau$. We also define $t_k=k\tau$. $\kappa>0$ is a parameter, which is be supposed to be small enough. We also control $\tau$: for some $\tau_0>0$, $\tau\leq \tau_0$.

\subsection{Galerkin approximation}
The first step of the proof is to consider finite dimensional approximations of the $H$-valued processes $(Y(t))_{t\in\R^+}$ and $(Y_k)_{k\in \N}$: if we fix $N\geq 1$, we define $(Y^{(N)}(t))_{t\in\R^+}$ and  $(Y_{k}^{(N)})_{k\in\N}$ by the equations
$$dY^{(N)}(t)=BY^{(N)}(t)dt+G_N(Y^{(N)}(t))dt+dW^{(N)}(t)$$
and
$$Y_{k+1}^{(N)}=Y_{k}^{(N)}+\tau BY_{k+1}^{(N)}+\tau G(Y_{k}^{(N)})+\sqrt{\tau}P_N\chi_{k+1},$$
with the initial conditions $Y_{t=0}^{(N)}=Y_{k=0}^{(N)}=P_Ny.$

The projection $P_N$ and the nonlinear coefficient $G_N$ have been defined above. $W^{(N)}=P_NW$ is a $N$-dimensional Wiener process on the subspace $H_N$. We remark that the above equations are well-defined on $H_N$ - which is a stable subspace of $B$.

The important and not difficult to prove result is the following: for any fixed $t\in \R^+$ and $k\in \N$, when $N\rightarrow +\infty$ we have
$$\E|Y(t)-Y^{(N)}(t)|^2\rightarrow 0 \quad \text{and} \quad\E|Y_k-Y_{k}^{(N)}|^2\rightarrow 0.$$

We need test functions $\Phi_N$ adapted to the finite-dimensional approximation: for any $N\geq 1$, by restriction we define $\Phi_N(y)=\Phi(y)$ for any $y\in H_N$; we obtain the following decomposition
\begin{align*}
\E\Phi(Y(m\tau))-\E\Phi(Y_m)&=\E\Phi(Y(m\tau))-\E\Phi(Y^{(N)}(m\tau))\\
&+\E\Phi_N(Y^{(N)}(m\tau))-\E\Phi_N(Y_{m}^{(N)})\\
&+\E\Phi(Y_{m}^{(N)})-\E\Phi(Y_m);
\end{align*}
the first and the third terms converge to $0$ when $N\rightarrow +\infty$. In the sequel, we prove an estimate of the second term, which is uniform with respect to dimension $N$; letting $N\rightarrow +\infty$ then yields an estimate on the left hand side.

Hence we work with the finite dimensional approximation, but we omit the parameter $N$. The constants appearing below are independent of $N$.

In Section \ref{sectusefulestim}, we prove some estimates on $Y(t)$ and $Y_m$, and in Section \ref{sectasymptu} we focus on the asymptotic behaviour of the processes.

\subsection{Some useful estimates}\label{sectusefulestim}
Bounds on moments of $Y_t$ and $Y_k$ can be proved, uniformly with respect to time.

\begin{lemme}\label{lem1}
For any $p\geq 1$, there exists a constant $C_p>0$ such that for every $t\geq 0$ and $y\in H$
$$\E|Y(t,y)|^p\leq C_p(1+|y|^p).$$
\end{lemme}
\underline{Proof}
If we define $Z(t)=Y(t)-W^B(t)$, we have $Z(0)=Y(0)=y$, and
$$\frac{dZ(t)}{dt}=BZ(t)+G(Y(t)),$$
and by Proposition \ref{propodiss}
\begin{align*}
\frac{1}{2}\frac{d|Z(t)|^2}{dt}&=<BZ(t)+G(Y(t)),Z(t)>\\
&=<BZ(t)+G(Z(t)),Z(t)>+<G(Y(t))-G(Z(t)),Z(t)>\\
&\leq -c|Z(t)|^2+C+\|G\|_{\infty}|Z(t)|\\
&\leq -c'|Z(t)|^2+C',
\end{align*}
for some new constants $c',C'$.

Then almost surely we have for any $t\geq 0$
$$|Z(t)|\leq C(1+|y|).$$
Thanks to \eqref{momWBp}, the conclusion easily follows.
\qed

\begin{lemme}\label{lem2}
For any $p\geq 1$, $\tau_0>0$, there exists a constant $C>0$ such that for every $0<\tau\leq \tau_0$, $k\in \N$ and $y\in H$
$$\E|Y_k|^{p}\leq C(1+|y|^p).$$
\end{lemme}

\underline{Proof}
As in the proof of Lemma \ref{lem1} above, we introduce $Z_m=Y_m-w_m$, where the process $(w_m)$ is the numerical approximation of $W^B$ with the numerical scheme \eqref{defYk} - with $G=0$; it is defined by
$$w_{m+1}=R_\tau w_m+\sqrt{\tau}R_\tau\chi_{m+1}.$$
Using Theorem $3.2$ of \cite{prin2}, giving the strong order $1/4$ for the numerical scheme - when the initial condition is $0$, with no nonlinear coefficient, with a constant diffusion term and under the assumptions made here - we obtain the following estimate: for any $p\geq 1$, $\tau_0>0$ and $0<r<1/2$ there exists $C>0$ such that for any $0<\tau\leq \tau_0$ and $m\geq 0$
\begin{equation}\label{estimnoisenum}
\E|w_{m}-W^B(m\tau)|^{2p}\leq C\tau^{(1/2-r)p}.
\end{equation}
Thanks to \eqref{momWBp} and \eqref{estimnoisenum}, we get that for any $\tau_0>0$, there exists $C>0$ such that for $0<\tau\leq \tau_0$ and $m\geq 0$
\begin{equation}\label{estimunifwm}
\E|w_{m}|^{2}\leq C.
\end{equation}
Now $Z_m$ defined above satisfies $Z_0=Y_0=y$ and
$$Z_{m+1}=R_\tau Z_m+\tau R_\tau G(Y_m);$$
since $|R_\tau|_{\mathcal{L}(H)}\leq \frac{1}{1+\mu_0\tau}$, we obtain the almost sure estimates
$$|Z_{m+1}|\leq \frac{1}{1+\mu_0\tau}|Z_m|+C\tau$$
and
$$|Z_m|\leq C(1+|y|).$$
Thanks to \eqref{estimunifwm}, we therefore obtain the result.
\qed

We now introduce the following process: for $0\leq k\leq m-1$ and $t_{k}\leq t\leq t_{k+1}$
\begin{equation}\label{tild}
\tilde{Y}(t)=Y_k+\int_{t_k}^{t}[B_\tau Y_k+R_\tau G(Y_k)]ds+\int_{t_k}^{t}R_\tau dW(s),
\end{equation}
where $B_\tau=BR_\tau$. The process $\tilde{Y}$ is a natural interpolation of the numerical solution $(Y_k)$ defined by \eqref{defYk}: $\tilde{Y}(t_k)=Y_k$.

Thanks to Lemma \ref{lem2}, we get
\begin{lemme}\label{lem2bis}
For any $p\geq 1$, $\tau_0>0$, there exists $C>0$ such that for any $0<\tau\leq \tau_0$, $t\geq 0$ and $y\in H$
$$\E|\tilde{Y}(t)|^p\leq C(1+|y|^p).$$
\end{lemme}

In the next Lemma, we give a control on Malliavin derivatives of $Y_k$ used in the proof. For an introduction to Malliavin calculus, see \cite{Nua}, \cite{Sanz}. The notations here are the same as in \cite{deb}, where the following useful integration by parts formula is given - see Lemma $2.1$ therein:
\begin{lemme}\label{lemintbyparts}
For any $F\in \mathbb{D}^{1,2}(H)$, $u\in \mathcal{C}_{b}^{2}(H)$ and $\Psi\in L^2(\Omega\times[0,T],\mathcal{L}_2(H))$ an adapted process,
\begin{equation}\label{formulaintbyparts}
\E[Du(F).\int_{0}^{T}\Psi(s)dW(s)]=\E[\int_{0}^{T}\text{Tr}(\Psi(s)^*D^2u(F)D_sF)ds],
\end{equation}
where $D_sF:h\in H\mapsto D_{s}^{h}F\in H$ stands for the Malliavin derivative of $F$, and $\mathbb{D}^{1,2}(H)$ is the set of $H$-valued random variables $F=\sum_{i\in\N} F_if_i$, with $F_i\in \mathbb{D}^{1,2}$ the domain of the Malliavin derivative for $\R$-valued random variables for any $i$, and $\sum_{i\in\N}\int_{0}^{T}\E|D_sF_i|^2ds<+\infty$.
\end{lemme}

Without any further dissipativity assumption, we are not able to give a uniform control with respect to time of the Malliavin derivative of $\tilde{Y}$. In the proof below, we take care of this problem by using this derivatives only at times $t_k=k\tau$ and $s$ such that $t_{k-l_s}\leq 1$.

\begin{lemme}\label{lem5}
For any $0\leq \beta<1$ and $\tau_0>0$, there exists a constant $C>0$ such that for every $h\in H$, $k\geq 1$, $0<\tau\leq \tau_0$ and $s\in \left[0,t_k\right]$
$$|D_{s}^{h}Y_k|_\beta\leq C(1+L_G\tau)^{k-l_s}(1+\frac{1}{(1+\mu_0\tau)^{(1-\beta)(k-l_s)}t_{k-l_s}^{\beta}})|h|.$$
\end{lemme}

\underline{Proof}
For any $k\geq 1$, $h\in H$ and $s\in\left[0,t_k\right]$, using the chain rule for Malliavin calculus and expressions \eqref{exprYk} and \eqref{exprsto}, we have
$$D_{s}^{h}Y_k=R_{\tau}^{k-l_s}h+\tau\sum_{l=l_s+1}^{k-1}R_{\tau}^{k-l}DG(Y_l).D_{s}^{h}Y_l.$$
Indeed, recall that $l_s$ denotes the integer part of $\frac{s}{\tau}$, so that when $l\leq l_s$ we have $D_{s}^{h}Y_l=0$.

As a consequence, for $k\geq l_s+1$
$$|D_{s}^{h}Y_k|\leq (1+L_G\tau)^{k-l_s}|h|.$$
Now using Lemma \ref{lem6}, we have
\begin{align*}
|(-B)^\beta D_{s}^{h}Y_k|&\leq \frac{1}{(1+\mu_0\tau)^{(1-\beta)(k-l_s)}t_{k-l_s}^{\beta}}|h|\\
&+L_G\tau\sum_{l=l_s+1}^{k-1}\frac{(1+L_G\tau)^{l-l_s}}{(1+\mu_0\tau)^{(1-\beta)(k-l)}t_{k-l}^{\beta}}|h|.
\end{align*}
To conclude, we see that
\begin{align*}
\tau\sum_{l=l_s+1}^{k-1}\frac{1}{(1+\mu_0\tau)^{(1-\beta)(k-l)}t_{k-l}^{\beta}}&\leq C\int_{0}^{+\infty}t^{-\beta}\frac{1}{(1+\mu_0\tau)^{(1-\beta)\frac{t}{\tau}}}dt\\
&\leq C<+\infty,
\end{align*}
when $0<\tau\leq \tau_0$.
\qed

\subsection{Asymptotic behaviour of the processes}\label{sectasymptu}

The results of this section are obtained for the initial $H$-valued processes, and for their finite dimensional approximations.

First, we focus on the existence of invariant measures for the continuous and discrete time processes. We use the well-known Krylov-Bogoliubov criterion - see \cite{DaP-Z2}. Tightness comes from two facts: $D(-B)^\gamma$ is compactly embedded in $H$ when $\gamma>0$, and when $\gamma<1/4$ we can control moments:
\begin{lemme}\label{estimtight}
For any $0<\gamma<1/4$, $\tau>0$ and any $y\in H$, there exists $C(\gamma,\tau,y),C(\gamma,y)>0$ such that for any $m\geq 1$ and $t\geq 1$
$$\E|Y_{m}(\tau,y)|_{\gamma}^{2}\leq C(\gamma,\tau,y) \quad \text{ and }\quad \E|Y(t,y)|_{\gamma}^{2}\leq C(\gamma,y)$$
\end{lemme}
Uniqueness of the invariant probability measure for the continuous time process $(Y(t))_{t\in\R^+}$ can be deduced from the well-known Doob Theorem - see \cite{DaP-Z2}. Indeed, since in equation \eqref{eqY} noise is additive and non-degenerate, the Strong Feller property - see also Lemma \ref{lem00} below - and irreducibility can be easily proved. In the proof of the main Theorem \ref{th}, we also need speed of convergence, and thanks to a coupling argument we get the following exponential convergence result - for a proof see Section $6.1$ in \cite{deb-hu-tess}:
\begin{propo}\label{propoexpy1y2}
There exist $c>0$, $C>0$ such that for any bounded test function $\phi$, any $t\geq 0$ and any $y_1,y_2\in H$
\begin{equation}\label{cvexpy1y2}
|\E\phi(Y(t,y_1))-\E\phi(Y(t,y_2))|\leq C\|\phi\|_{\infty}(1+|y_1|^2+|y_2|^2)e^{-ct}.
\end{equation}
\end{propo}

The idea of coupling relies on the following formula: if $\nu_1$ and $\nu_2$ are two probability measures on a state space $S$, their total variation distance satisfies
$$d_{TV}(\nu_1,\nu_2)=\inf\left\{\PP(X_1\neq X_2)\right\},$$
which is an infimum over random variables $(X_1,X_2)$ defined on a same probability space, and such that $X_1\sim\nu_1$ and $X_2\sim \nu_2$.

Roughly speaking, the principle is to define a coupling $(Z_1(t,y_1,y_2),Z_2(t,y_1,y_2))_{t\geq 0}$ for the processes $(Y(t,y_1)_{t\geq 0}$ and $Y((t,y_2))_{t\geq 0}$ such that the coupling time $\mathcal{T}$ of $Z_1$ and $Z_2$ - i.e. the first time the processes are equal - has an exponentially decreasing tail.

This technique was first used in the study of the asymptotic behaviour of Markov chains - see \cite{Brem}, \cite{Doe}, \cite{Lin}, \cite{Me-Twee} - and was later adapted for SDEs and more recently for SPDEs - see for instance \cite{ku-shi}, \cite{matt}.


In fact, uniqueness of an invariant probability measure $\overline{\mu}$ is an easy consequence of this Proposition, and moreover we get for any $y\in H$ and any $t\geq 0$
\begin{equation}\label{cvexpy_int}
|\E\phi(Y(t,y))-\int_{H}\phi d\overline{\mu}|\leq C\|\phi\|_{\infty}(1+|y|^2)e^{-ct}.
\end{equation}

In general, we do not know whether uniqueness also holds for the numerical approximation $(Y_k(\tau,.))_{k\in \N}$.


\begin{rem}\label{remstrictergo}
A sufficient condition for the uniqueness of the invariant probability measure of the discrete time process $(Y_k)_{k\in \N}$ is the strict dissipativity assumption
$$L_G<\mu_0,$$
where we recall that $L_G$ denotes the Lipschitz constant of $G$.

Then trajectories of the processes $(Y_t)_{t\in\R^+}$ and $(Y_k)_{k\in\N}$ issued from different initial conditions $y_1$ and $y_2$ and driven by the same noise process are exponentially close when time increases: for any $\tau_0>0$, there exists $c>0$ such that for any $0<\tau\leq \tau_0$, $k\geq 0$ and $t\geq 0$ we have almost surely
\begin{gather*}
|Y(t,y_1)-Y(t,y_2)|\leq e^{-(\mu_0-L_G)t}|y_1-y_2|\\
|Y_k(\tau,y_1)-Y_k(\tau,y_2)|\leq e^{-ck\tau}|y_1-y_2|.
\end{gather*}
Proof of uniqueness is then straightforward - and we do not need Proposition \ref{propoexpy1y2} above.
\end{rem}

\section{Presentation of the proof of the weak approximation result}\label{sectproof}

The proof of Theorem \ref{th} is very technical, so for pedagogy we first introduce the decomposition of the error, and identify the term which we control later in Section \ref{sectestims}. Some crucial estimates on the derivatives of the semi-group with respect to the initial conditions - regularization, long-time behaviour - are proved below in Sub-Section \ref{expdec}.

\subsection{Strategy}

We define
\begin{equation}\label{defu}
u(t,y)=\E[\phi(Y(t,y))],
\end{equation}
which is solution of a finite dimensional Kolmogorov equation associated with the finite dimensional approximation of \eqref{eqY}:
$$\frac{\text{d}u}{\text{d}t}(t,y)=Lu(t,y)=\frac{1}{2}\text{Tr}\left(D^2u(t,y)\right)+<By+G(y),Du(t,y)>.$$
As explained in the introduction, this is one of the essential tools in the proof of the weak approximation result.

The weak error at time $T=m\tau$ can be decomposed with a telescoping sum - where to simplify the dependence of the numerical approximation in $\tau$ and $y$ is not written:
\begin{equation}\label{decompabc_k}
\begin{aligned}
\E[\phi(Y(T,y))]-\E[\phi(Y_m(\tau,y))]&=u(T,y)-\E[u(0,Y_m(\tau,y))]\\
&=\sum_{k=0}^{m-1}(\E[u(T-t_k,Y_k)]-\E[u(T-t_{k+1},Y_{k+1})])\\
&=u(T,y)-\E[u(T-\tau,Y_1(\tau,y))]+\sum_{k=1}^{m-1}(a_k+b_k+c_k),
\end{aligned}
\end{equation}
where for $1\leq k\leq m-1$
\begin{equation}\label{defabc_k}
\begin{gathered}
a_k=\E\int_{t_k}^{t_{k+1}}<B\tilde{Y}(t)-B_\tau Y_k,Du(T-t,\tilde{Y}(t))>dt,\\
b_k=\E\int_{t_k}^{t_{k+1}}<G(\tilde{Y}(t))-R_\tau G(Y_k),Du(T-t,\tilde{Y}(t))>dt,\\
c_k=\frac{1}{2}\E\int_{t_k}^{t_{k+1}}\text{Tr}((I-R_\tau R_{\tau}^{*})D^2u(T-t,\tilde{Y}(t)))dt.
\end{gathered}
\end{equation}
This follows from the use of the Kolmogorov equation and of the It\^{o} formula.

\subsection{Bounds on the derivatives of the transition semi-group}\label{expdec}

By \eqref{defu}, $u(t,y)=\E[\phi(Y(t,y))]$; since $\phi$ is of class $\mathcal{C}^2$, bounded and with bounded derivatives, we are able to prove that with respect to $y$ the function $u$ is twice differentiable, and that the derivatives can be calculated in the following way:
\begin{itemize}
\item For any $h\in H$, we have
\begin{equation}\label{formDu}
Du(t,y).h=\E[D\phi(Y(t,y)).\eta^{h,y}(t)],
\end{equation}
where $\eta^{h,y}$ is the solution of
\begin{gather*}
\frac{d\eta^{h,y}(t)}{dt}=B\eta^{h,y}(t)+DG(Y(t,y)).\eta^{h,y}(t),\\
\eta^{h,y}(0)=h.
\end{gather*}
\item For any $h,k\in H$, we have
\begin{equation}\label{formD^2u}
D^2u(t,y).(h,k)=\E[D^2\phi(Y(t,y)).(\eta^{h,y}(t),\eta^{k,y}(t))+D\phi(Y(t,y)).\zeta^{h,k,y}(t)],
\end{equation}
where $\zeta^{h,k,y}$ is the solution of
\begin{gather*}
\frac{d\zeta^{h,k,y}(t)}{dt}=B\zeta^{h,k,y}(t)+DG(Y(t,y)).\zeta^{h,k,y}(t)+D^2G(Y(t,y)).(\eta^{h,y}(t),\eta^{k,y}(t)),\\
\zeta^{h,k,y}(0)=0.
\end{gather*}
\end{itemize}

In \cite{deb}, the key point for obtaining the weak order $1/2$ is to control the derivatives $|Du(t,y)|_{\beta}$ and $|(-B)^\beta D^2 u(t,y)(-B)^\gamma|_{\mathcal{L}(H)}$, with $\beta<1/2$ and $\gamma<1/2$ - with the identification of Remark \ref{remident}. Moreover, to obtain a long-time weak estimate we need to prove some exponential decreasing of such quantities when time $t$ goes to infinity. The two Propositions below are the essential results we thus need.

\begin{propo}\label{lem3}
There exists a constant $\tilde{\mu}>0$ such that for any $0\leq\beta<1/2$, for any $t>0$ and $y\in H$
\begin{equation}\label{decreaseDU2}
|Du(t,y)|_{\beta}\leq C_\beta(1+\frac{1}{t^\beta})e^{-\tilde{\mu}t}(1+|y|^2).
\end{equation}
\end{propo}

\begin{propo}\label{lem4}
There exists a constant $\tilde{\mu}>0$ such that for any $0\leq\beta,\gamma<1/2$, for any $t>0$ and $y\in H$
\begin{equation}\label{decreaseDU3}
|(-B)^\beta D^2 u(t,y)(-B)^\gamma|_{\mathcal{L}(H)}\leq C_{\beta,\gamma} (1+\frac{1}{t^\eta}+\frac{1}{t^{\beta+\gamma}})e^{-\tilde{\mu}t}(1+|y|^2).
\end{equation}
\end{propo}

The singularity $t^{-\eta}$ in \eqref{decreaseDU3} is a consequence of the regularity properties satisfied by $G$. Since in general during the proof of Theorem \ref{th}, we need $\beta+\gamma$ to be close to $1$, and therefore greater than $\eta$, only the second singularity $t^{-\beta-\gamma}$ plays a role.

The proofs require several steps. First in Lemma \ref{lem000} below we prove estimates for a finite horizon and general $0\leq\beta,\gamma<1/2$; then in Lemma \ref{lem00} we study the long-time behaviour in the particular case $\beta=\gamma=0$; we finally conclude with the proofs of Propositions \ref{lem3} and \ref{lem4}.

First, we prove estimates of these quantities for $0<t\leq 1$ - see Lemmas $4.4$ and $4.5$ in \cite{deb}, with a difference coming from the assumptions made on the nonlinear coefficient $G$:
\begin{lemme}\label{lem000}
For any $0\leq \beta<1/2$, $0\leq \gamma<1/2$, there exist constants $C_{\beta}$ and $C_{\beta,\gamma}$ such that for any $y\in H$ and any $0<t\leq 1$
\begin{gather*}
|Du(t,y)|_{\beta}\leq \frac{C_\beta}{t^\beta}\\
|(-B)^\beta D^2 u(t,y)(-B)^\gamma|_{\mathcal{L}(H)}\leq C_{\beta,\gamma}(\frac{1}{t^\eta}+\frac{1}{t^{\beta+\gamma}}).
\end{gather*}
\end{lemme}
If we take another time interval $]0,T_{max}]$ instead of $]0,1]$, the constants $C_{\beta}$ and $C_{\beta,\gamma}$ are \textit{a priori} exponentially increasing in $T_{max}$.

\underline{Proof}
Owing to \eqref{formDu} and \eqref{formD^2u}, we only need to prove the following almost sure estimates, for some constants $C_\beta$ and $C_{\beta,\gamma}$ - which may vary from line to line below: for any $0<t\leq 1$
\begin{equation}\label{deriv_imp}
\begin{gathered}
|\eta^{h,y}(t)|\leq \frac{C_\beta}{t^\beta}|h|_{-\beta}\\
|\zeta^{h,k,y}(t)|\leq \frac{C_{\beta,\gamma}}{t^\eta}|h|_{-\beta}|k|_{-\gamma},
\end{gathered}
\end{equation}
where the parameter $\eta$ is defined in Assumption \ref{hypG}.

We use mild formulations, and the regularization properties of the semi-group given in Proposition \ref{proporegul}:
\begin{align*}
|\eta^{h,y}(t)|&=|e^{tB}h+\int_{0}^{t}e^{(t-s)B}DG(Y(s,y)).\eta^{h,y}(s)ds|\\
&\leq \frac{C_\beta}{t^\beta}|h|_{-\beta}+C\int_{0}^{t}|\eta^{h,y}(s)|ds,
\end{align*}
and by the Gronwall Lemma we get the result.

For the second-order derivative, we moreover use the properties of $G$ in Assumption \ref{hypG} to get
\begin{align*}
|\zeta^{h,k,y}(t)|&=|\int_{0}^{t}e^{(t-s)B}DG(Y(s,y)).\zeta^{h,k,y}(s)ds\\
&+\int_{0}^{t}e^{(t-s)B}D^2G(Y(s,y)).(\eta^{h,y}(s),\eta^{k,y}(s))ds|\\
&\leq C\int_{0}^{t}|\zeta^{h,k,y}(s)|ds+\int_{0}^{t}\frac{C_{\beta,\gamma}}{(t-s)^\eta}|\eta^{h,y}(s)||\eta^{k,y}(s)|ds\\
&\leq C\int_{0}^{t}|\zeta^{h,k,y}(s)|ds+C_{\beta,\gamma}|h|_{-\beta}|k|_{-\beta}t^{1-\eta-\beta-\gamma}\int_{0}^{1}\frac{1}{(1-s)^\eta s^{\beta+\gamma}}ds.
\end{align*}
To conclude, it remains to use the Gronwall Lemma.
\qed

Thanks to the dissipativity property expressed in Proposition \ref{propodiss}, we can prove the result in the case $\beta=\gamma=0$. We notice that the proof would be straightforward under a strict dissipativity assumption - since then $\eta^{h,y}(t)$ and $\zeta^{h,k,y}(t)$ would decrease exponentially in $t$; in the general case $\eta^{h,y}(t)$ and $\zeta^{h,k,y}(t)$ are exponentially increasing in time so that we can not work directly. Here the result comes from the estimate \eqref{cvexpy1y2} of Proposition \ref{propoexpy1y2}.
\begin{lemme}\label{lem00}
There exist constants $C$ and $c>0$ such that for any $t\geq 0$ and any $y\in H$
\begin{equation}\label{decreaseDU1}
|Du(t,y)|\leq Ce^{-ct}(1+|y|^2) \quad \text{ and } \quad |D^2u(t,y)|_{\mathcal{L}(H)}\leq Ce^{-ct}(1+\frac{1}{t^\eta})(1+|y|^2).
\end{equation}
\end{lemme}

\underline{Proof}
The Bismut-Elworthy-Li formula states that if $\Phi:H\rightarrow \R$ is a function of class $\mathcal{C}^2$ with bounded derivatives and with at most quadratic growth - i.e. there exists $M(\Phi)>0$ such that for any $y\in H$ we have $|\Phi(y)|\leq M(\Phi)(1+|y|^2)$ - then we can calculate the first and the second order derivatives of $(t,y)\mapsto v(t,y):=\E\Phi(Y(t,y))$ with respect to $y$.
First, we have for any $y\in H$ and $h\in H$
\begin{equation}\label{BEL1}
\begin{aligned}
Dv(t,y).h&=\frac{1}{t}\E[\int_{0}^{t}<\eta^{h,y}(s),dW(s)>\Phi(Y(t,y))]\\
&=\frac{2}{t}\E[\int_{0}^{t/2}<\eta^{h,y}(s),dW(s)>v(t/2,Y(t/2,y))];
\end{aligned}
\end{equation}
the second equality is a consequence of the identity $v(t,y)=\E v(t/2,Y(t/2,y))$, thanks to the Markov property.
Using the second formula of \eqref{BEL1}, we obtain a formula for the second order derivative: for any $y\in H$ and $h,k\in H$,
\begin{equation}\label{BEL2}
\begin{aligned}
D^2v(t,y).(h,k)&=\frac{2}{t}\E[\int_{0}^{t/2}<\zeta^{h,k,y}(s),dW(s)>v(t/2,Y(T/2,y))]\\
&+\frac{2}{t}\E[\int_{0}^{t/2}<\eta^{h,y}(s),dW(s)>Dv(t/2,Y(t/2)).\eta^{k,y}(t/2)].
\end{aligned}
\end{equation}
We then see, using Lemmas \ref{lem1} and \ref{lem000} - with $\beta=\gamma=0$ - that there exists $C>0$ such that for any $0<t\leq 1$, $y\in H$, $h,k\in H$
\begin{equation}\label{consBEL}
\begin{gathered}
|Dv(t,y).h|\leq \frac{C}{\sqrt{t}}M(\Phi)(1+|y|^2)|h|,\\
|D^2v(t,y).(h,k)|\leq \frac{C}{t}M(\Phi)(1+|y|^2)|h||k|.
\end{gathered}
\end{equation}
Now when $t\geq 1$ the Markov property implies that $u(t,y)=\E u(t-1,Y(1,y))$, and by \eqref{cvexpy_int} we have
$$|u(t-1,y)-\int_H\phi d\mu|\leq Ce^{-c(t-1)}(1+|y|^2).$$
If we choose $\Phi_t(y)=u(t-1,y)-\int_H\phi d\mu$, we have $u(t,y)=\E\Phi_t(Y(1,y))+\int_H\phi d\mu$, with $M(\Phi_t)\leq Ce^{-c(t-1)}$. With \eqref{consBEL} at time $1$, we obtain for $t\geq 1$
\begin{gather*}
|Du(t,y).h|\leq Ce^{-c(t-1)}(1+|y|^2)|h|\\
|D^2u(t,y).(h,k)|\leq Ce^{-c(t-1)}(1+|y|^2)|h||k|.
\end{gather*}
Moreover by Lemma \ref{lem000} we have a control when $0\leq t\leq 1$, so that with a change of constants we get the result.
\qed



We can finally prove the Propositions \ref{lem3} and \ref{lem4}. The key tool is the Markov property of the process $Y$ which yields the following formula: for any $t\geq 1$
\begin{equation}\label{umarkov}
u(t,y)=\E[u(t-1,Y_1(y))].
\end{equation}

To get the exponential decreasing, we use Lemma \ref{lem00} at time $t-1$ when $t\geq 1$, while $|h|_{-\beta}$ appears from $\eta_{h,y}(1)$, and with estimates coming from Lemma \ref{lem000}.

\vspace{0.5cm}

\underline{Proof of Propositions \ref{lem3} and \ref{lem4}}
Using equation \eqref{umarkov} and Lemma \ref{lem00}, for any $t\geq 1$ we have
$$|Du(t,y).h|\leq Ce^{-c(t-1)}\E[(1+|Y(1,Y)|^2)|\eta^{h,y}(1)|]\leq Ce^{-c(t-1)}(1+|y|^2)|h|_{-\beta},$$
where the last estimate comes from Lemmas \ref{lem1} and \ref{lem000}.
Combining this estimate with the result of Lemma \ref{lem000}, which gives an estimate for $t\leq 1$, we obtain \eqref{decreaseDU2}.
For the second order derivatives, Lemma \ref{lem000} gives an estimate for $t\leq 1$, and for $t\geq 1$ we use \eqref{umarkov} to see that
\begin{align*}
D^2u(t,y).(h,k)&=\E[D^2[u(t-1,Y(1,y))].(h,k)]\\
&=\E D^2u(t-1,Y(1,y)).(\eta^{h,y}(1),\eta^{k,y}(1))+\E Du(t-1,Y(1,y)).\zeta^{h,k,y}(1).
\end{align*}
Using Lemma \ref{lem00}, we get an exponential decreasing; thanks to Lemma \ref{lem1} and to the estimates in the proof of Lemma \ref{lem000} at time $1$, we obtain
$$|D^2u(t,y).(h,k)|\leq Ce^{-c(t-1)}(1+|y|^2)|h|_{-\beta}|k|_{-\gamma}.$$
Then \eqref{decreaseDU3} easily follows.
\qed


\subsection{Proof of Corollary \ref{corMesInv}}


The first estimate is a simple consequence of the Theorem, and of the exponential convergence to equilibrium of the continuous-time process - see \eqref{cvexpy_int}. We then get

$$|\E[\phi(Y_{m}(\tau,y))]-\int_{H}\phi d\overline{\mu}|\leq C(1+|y|^3)(\frac{1}{m^{1/2-\kappa}}+\tau^{1/2-\kappa})+C(1+|y|^2)e^{-cm\tau}.$$



If $\mu^\tau$ is an ergodic invariant probability measure of $(Y_m(\tau,.))_m$, then since $\phi$ is bounded for $\mu^\tau$-almost any $y\in H$ we have by the ergodic Theorem the following convergence when $M\rightarrow +\infty$:
$$\frac{1}{M}\sum_{m=1}^{M}\E[\phi(Y_{m}(\tau,y))]\rightarrow \int_H\phi(y)\mu^\tau(dy).$$
To conclude, it remains to choose a initial condition $y$ in this non-empty set, and to use Cesaro Lemma on the right-hand side of the estimate.

We notice that if $\mu^\tau$ is an invariant probability measure, not necessarily ergodic, having a finite moment of order $3$, then it is enough to integrate the inequality above with respect to $\mu^\tau$.


\section{Proof of the estimates}\label{sectestims}

We need to control the terms given in \eqref{defabc_k}, according to the decomposition \eqref{decompabc_k}. We recall that constants $C$ must be independent from the dimension $N$ and the final time $T=m\tau$.

\subsection{Estimate of $u(T,x)-\E[u(T-\tau,Y_1)]$}\label{lackreg}


The Markov property gives
$$u(T,y)=\E[\phi(Y(T,y))]=\E[u(T-\tau,Y(\tau,y))].$$
If $0<\kappa<1/2$, using Lemma \ref{lem2} and Proposition \ref{lem3} we get
$$|u(T,y)-\E[u(T-\tau,Y_1)]|\leq C(1+(T-\tau)^{-1/2+\kappa})e^{-\tilde{\mu}(T-\tau)}(\E|Y(\tau,y)-Y_1|_{-1/2+\kappa}^{2})^{1/2}(1+|y|^2).$$
We can write
\begin{align*}
Y(\tau,y)-Y_1&=(e^{\tau B}-R_\tau)y+\int_{0}^{\tau}e^{(\tau-s)B}G(Y(s,y))ds-\tau R_\tau G(y)\\
&+\int_{0}^{\tau}e^{(\tau-s)B}dW(s)-\sqrt{\tau}R_\tau\chi_1.
\end{align*}
We use the following properties to estimate the first line in this equality:
\begin{gather*}
|(-B)^{-1/2+\kappa}(e^{\tau B}-R_\tau)|_{\mathcal{L}(H)}\leq c\tau^{1/2-\kappa},\\
|e^{sB}|_{\mathcal{L}(H)}\leq 1 \text{ for }s\geq 0,\\
|R_\tau|_{\mathcal{L}(H)}\leq 1,\\
|(-B)^{-1/2+\kappa}.|\leq c|.|;
\end{gather*}
therefore the first line in the last expression is almost surely bounded by $C(\tau^{1/2-\kappa}+\tau)(1+|y|).$

For the second line, we have
\begin{align*}
\E|(-B)^{-1/2+\kappa}\int_{0}^{\tau}e^{(\tau-s)B}dW(s)|^2&=\E\int_{0}^{\tau}|(-B)^{-1/2+\kappa}e^{(\tau-s)B}|_{\mathcal{L}_2(H)}^2ds\\
&\leq \tau |(-B)^{-1/2+\kappa}|_{\mathcal{L}_2(H)}^2\\
&\leq c\tau;
\end{align*}
the last term is controlled in the same way: $\E|(-B)^{-1/2+\kappa}\sqrt{\tau}R_\tau\chi_1|^2\leq c\tau$. Therefore we have
\begin{equation}\label{e1}
|u(T,x)-\E[u(T-\tau,Y_1)]|\leq C(1+|y|^3)(1+(T-\tau)^{-1/2+\kappa})e^{-\tilde{\mu}(T-\tau)}\tau^{1/2-\kappa}.
\end{equation}

We thus understand that to obtain weak order $1/2$ requires to be careful in the estimate. Here we used Lemma \ref{lem2} instead of Lemma \ref{lem00}; otherwise looking at $\E|Y(\tau,y)-Y_1|^2$ would have not been sufficient. The control of the other terms must be done in the same spirit.

\subsection{Estimate of $a_k$}

We have
\begin{align*}
a_k&=\E\int_{t_k}^{t_{k+1}}<B\tilde{Y}(t)-B_\tau Y_k,Du(T-t,\tilde{Y}(t))>dt\\
&=\E\int_{t_k}^{t_{k+1}}<(B-B_\tau)Y_k,Du(T-t,\tilde{Y}(t))>dt\\
&+\E\int_{t_k}^{t_{k+1}}<B(\tilde{Y}(t)-Y_k),Du(T-t,\tilde{Y}(t))>dt\\
&:=a_{k}^{1}+a_{k}^{2}.
\end{align*}

\subsubsection{Estimate of $a_{k}^{1}$}\label{pourint}

We use the equality $B_\tau-B=\tau R_\tau B^2$. We also decompose $a_{k}^{1}$ using expression \eqref{exprYk}:
\begin{gather*}
a_{k}^{1,1}=-\tau\E\int_{t_k}^{t_{k+1}}<R_\tau B^2R_{\tau}^{k}y,Du(T-t,\tilde{Y}(t))>dt\\
a_{k}^{1,2}=-\tau\E\int_{t_k}^{t_{k+1}}<R_\tau B^2\tau\sum_{l=0}^{k-1}R_{\tau}^{k-l}G(Y_l),Du(T-t,\tilde{Y}(t))>dt\\
a_{k}^{1,3}=-\tau\E\int_{t_k}^{t_{k+1}}<R_\tau B^2\sqrt{\tau}\sum_{l=0}^{k-1}R_{\tau}^{k-l}\chi_{l+1},Du(T-t,\tilde{Y}(t))>dt;
\end{gather*}
then $a_{k}^{1}=a_{k}^{1,1}+a_{k}^{1,2}+a_{k}^{1,3}$.

\begin{enumerate}
\item \textbf{Estimate of $a_{k}^{1,1}$}

The idea is to ``share'' $B^2$ between different factors - thanks to regularization properties of the semi-group $(R_{\tau}^{k})_{k\in \N}$ and to Lemma \ref{lem2} - in order to increase the order of convergence.

\begin{align*}
|a_{k}^{1,1}|&\leq \tau\E\int_{t_k}^{t_{k+1}}|R_\tau(-B)^{1/2+2\kappa}|_{\mathcal{L}(H)}|(-B)^{1-\kappa}R_{\tau}^{k}|_{\mathcal{L}(H)}|y|_H|(-B)^{1/2-\kappa}Du(T-t,\tilde{Y}(t))|dt\\
&\leq C|y|\tau \tau^{-1/2-2\kappa}t_{k}^{-1+\kappa}\int_{t_k}^{t_{k+1}}(1+(T-t)^{-1/2+\kappa})e^{-\tilde{\mu}(T-t)}\E(1+|\tilde{Y}(t)|^2)dt,
\end{align*}
thanks to Lemma \ref{lem6} and Proposition \ref{lem3}.

By taking expectation, thanks to Lemma \ref{lem2bis} we have
\begin{align*}
\sum_{k=1}^{m-1}|a_{k}^{1,1}|&\leq C|y|\tau^{1/2-2\kappa}\sum_{k=1}^{m-1}\frac{1}{t_{k}^{1-\kappa}}\int_{t_{k}}^{t_{k+1}}(1+\frac{1}{(T-t)^{1/2-\kappa}})e^{-\tilde{\mu}(T-t)}(1+|y|^2)dt\\
&\leq C_\kappa(1+|y|^3)\tau^{1/2-2\kappa}\int_{0}^{T}\frac{1}{t^{1-\kappa}(T-t)^{1/2-\kappa}}dt\\
&\leq T^{-(1/2-2\kappa)}C_\kappa(1+|y|^3)\tau^{1/2-2\kappa}\int_{0}^{1}\frac{1}{s^{1-\kappa}(1-s)^{1/2-\kappa}}ds,
\end{align*}
and we obtain
\begin{equation}\label{e2}
\sum_{k=1}^{m-1}|a_{k}^{1,1}|\leq C(1+|y|^3)T^{-(1/2-2\kappa)}\tau^{1/2-2\kappa}.
\end{equation}

\item \textbf{Estimate of $a_{k}^{1,2}$}

First we write
$$|a_{k}^{1,2}|\leq C\tau\E\int_{t_{k}}^{t_{k+1}}|R_\tau(-B)^{1/2+2\kappa}|_{\mathcal{L}(H)}|\tau(-B)^{1-\kappa}\sum_{l=0}^{k-1}R_{\tau}^{k-l}G(Y_l)||(-B)^{1/2-\kappa}Du(T-t,\tilde{Y}(t))|dt.$$

Using Lemma \ref{lem6}, we can prove the following useful inequality: for $\tau\leq \tau_0$ and any $k\geq 1$
\begin{equation}\label{sumestim}
\tau\sum_{l=1}^{k}\frac{1}{(l\tau)^{1-\kappa}}\frac{1}{(1+\mu_0\tau)^{l\kappa}}\leq C_\kappa.
\end{equation}

Indeed,
\begin{align*}
\tau\sum_{l=1}^{k}\frac{1}{(l\tau)^{1-\kappa}}\frac{1}{(1+\mu_0\tau)^{l\kappa}}&\leq C\int_{0}^{t_k}\frac{1}{t^{1-\kappa}}\frac{1}{(1+\mu_0\tau)^{\kappa\frac{t}{\tau}}}dt\\
&\leq \int_{0}^{\infty}\frac{1}{t^{1-\kappa}}e^{-t\frac{\kappa}{\tau}\log(1+\mu_0\tau)}dt\\
&\leq \int_{0}^{\infty}\frac{1}{s^{1-\kappa}}e^{-s}ds \left(\frac{\tau}{\kappa\log(1+\mu_0\tau)}\right)^\kappa\\
&\leq C_\kappa.
\end{align*}

Since $G$ is supposed to be bounded, the estimate \eqref{sumestim} yields
$$
|\tau(-B)^{1-\kappa}\sum_{l=0}^{k-1}R_{\tau}^{k-l}G(Y_l)|\leq C\|G\|_{\infty}\tau\sum_{l=1}^{k}\frac{1}{(l\tau)^{1-\kappa}}\frac{1}{(1+\mu_0\tau)^{l\kappa}}\leq C_\kappa.
$$

With Lemma \ref{lem2bis} and Proposition \ref{lem4}, we can now write
$$|a_{k}^{1,2}|\leq C
(1+|y|^2)\tau^{1/2-2\kappa}\int_{t_{k}}^{t_{k+1}}(1+\frac{1}{(T-t)^{1/2-\kappa}})e^{-\tilde{\mu}(T-t)}dt,$$
and we get
\begin{equation}\label{e3}
\sum_{k=1}^{m-1}|a_{k}^{1,2}|\leq C(1+|y|^2)\tau^{1/2-2\kappa}.
\end{equation}

\item \textbf{Estimate of $a_{k}^{1,3}$}

The analysis of this term is more complicated. We recall that since noise is white in space, for any $t>0$, $\E|(-B)^\gamma W^B(t)|^2<+\infty$ if and only if $\gamma<1/4$; as a consequence, the strategy used above to control $a_{k}^{1,1}$ can not be used directly - otherwise we could only obtain order $1/4$.

In \cite{deb}, an integration by parts formula is used to deal with the lack of regularity of the stochastic integral appearing in the definition of $a_{k}^{1,3}$. An additionnal difficulty arises in our situation because the estimate given in Lemma \ref{lem5} is not uniform with respect to time. Instead, we remark that the two problems - lack of regularity and bad time dependence - do not occur at the same time; therefore a decomposition of the interval $[0,t_k]$ into $[0,t_k-1]$ and $[t_k-1,t_k]$ - for $k$ large enough - can help to treat the problems separately.

Let us explain more concretely the situation at the continuous time level: we have to treat - at the finite dimensional approximation level, but with a bound independent from dimension - an expression involving $B^2\int_{0}^{t}e^{(t-s)B}dW(s)$. The idea to get rid of this expression is to use an integration by parts formula on the whole interval. We can also see that we can do this integration by parts only on a subinterval of size independent from $t$: indeed if $t\geq 1$
$$\int_{0}^{t}e^{(t-s)B}dW(s)=\int_{0}^{t-1}e^{(t-s)B}dW(s)+\int_{t-1}^{t}e^{(t-s)B}dW(s).$$
The first term is equal to $e^B\int_{0}^{t-1}e^{(t-1-s)B}dW(s)$, so that thanks to regularization properties of the semi-group $(e^{tB})$ - see Proposition \ref{proporegul} - multiplication by $B^2$ is possible - in other words we do not require an integration by parts to get a bound independent from the dimension; to treat the second term, the lack of regularity remains but can still be treated by the integration by parts, with the advantage of involving a smaller interval, where at the discrete time level below we can use a uniform control of all quantities, thanks to Lemma \ref{lem5}. Then the same idea of ``sharing'' $B^2$ can be used again to get order $1/2$.


Let us now explain develop this program for the discrete time situation: by using \eqref{exprsto}, we make the decomposition
\begin{align*}
a_{k}^{1,3}&=-\tau\E\int_{t_k}^{t_{k+1}}<R_\tau B^2\sqrt{\tau}\sum_{l=0}^{k-1}R_{\tau}^{k-l}\chi_{l+1},Du(T-t,\tilde{Y}(t))>dt\\
&=-\tau\E\int_{t_k}^{t_{k+1}}<\int_{0}^{t_k}R_\tau B^2R_{\tau}^{k-l_s}dW(s),Du(T-t,\tilde{Y}(t))>dt\\
&=-\tau\E\int_{t_k}^{t_{k+1}}<\int_{0}^{(t_k-1)\vee 0}R_\tau B^2R_{\tau}^{k-l_s}dW(s),Du(T-t,\tilde{Y}(t))>dt\\
&-\tau\E\int_{t_k}^{t_{k+1}}<\int_{(t_k-1)\vee 0}^{t_k}R_\tau B^2R_{\tau}^{k-l_s}dW(s),Du(T-t,\tilde{Y}(t))>dt.
\end{align*}

For the first term - which is equal to $0$ when $t_k<1$ - we use the Cauchy-Schwarz inequality and we directly get
\begin{align*}
|\E<\int_{0}^{(t_k-1)\vee 0}&R_\tau B^2R_{\tau}^{k-l_s}dW(s),Du(T-t,\tilde{Y}(t))>|\\
&\leq (\E|\int_{0}^{(t_k-1)\vee 0}R_\tau B^2R_{\tau}^{k-l_s}dW(s)|^2)^{1/2}(\E|Du(T-t,\tilde{Y}(t))|^2)^{1/2}\\
&\leq C(1+|y|^2)e^{-c(T-t)},
\end{align*}
thanks to Lemmas \ref{lem00} and \ref{lem2bis}, and to the following inequality - we remark that in the integral below $t_{k-l_s}\geq 1$:
\begin{align*}
\E|\int_{0}^{(t_k-1)\vee 0}R_\tau B^2R_{\tau}^{k-l_s}dW(s)|^2&=\int_{0}^{(t_k-1)\vee 0}|R_\tau B^2R_{\tau}^{k-l_s}|_{\mathcal{L}_2(H)}^{2}ds\\
&=\int_{0}^{(t_k-1)\vee 0}\text{Tr}(R_{\tau}^{2} B^4R_{\tau}^{2(k-l_s)})ds\\
&\leq \int_{0}^{(t_k-1)\vee 0}|(R_{\tau}^{2} B^{4+1/2+\kappa}R_{\tau}^{2(k-l_s)}|_{\mathcal{L}(H)}ds\text{Tr}((-B)^{-1/2-\kappa})\\
&\leq C\int_{0}^{(t_k-1)\vee 0}\frac{1}{(1+\mu_0\tau)^{k-l_s}t_{k-l_s}^{4+1/2+\kappa}}ds\\
&\leq C\int_{0}^{(t_k-1)\vee 0}\frac{1}{(1+\mu_0\tau)^{k-l_s}}ds\\
&\leq C\int_{0}^{+\infty}\frac{1}{(1+\mu_0\tau)^{s/\tau}}ds\\
&\leq C,
\end{align*}
when $\tau\leq \tau_0$.
Then
$$|\tau\E\int_{t_k}^{t_{k+1}}<\int_{0}^{(t_k-1)\vee 0}R_\tau B^2R_{\tau}^{k-l_s}dW(s),Du(T-t,\tilde{Y}(t))>dt|\leq C\int_{t_k}^{t_{k+1}}e^{-\tilde{\mu}(T-t)}dt(1+|y|^2)\tau.$$

For the second term, we use the integration by parts formula of Lemma \ref{lemintbyparts} to get
\begin{align*}
\tau\E\int_{t_k}^{t_{k+1}}&<\int_{(t_k-1)\vee 0}^{t_k}R_\tau B^2R_{\tau}^{k-l_s}dW(s),Du(T-t,\tilde{Y}(t))>dt\\
&=-\tau\E\int_{t_k}^{t_{k+1}}\int_{(t_k-1)\vee 0}^{t_k}\text{Tr}\left(R_{\tau}^{k-l_s}B^2R_\tau D^2u(T-t,\tilde{Y}(t))D_s\tilde{Y}(t)\right)dsdt.
\end{align*}
If $h\in H$, $(t_k-1)\vee 0\leq s\leq t_k\leq t< t_{k+1}$, with \eqref{tild} we see that
\begin{align*}
D_{s}^{h}\tilde{Y}(t)&=D_{s}^{h}Y_k+\int_{t_k}^{t}(B_\tau D_{s}^{h}Y_k+R_\tau G(Y_k)D_{s}^{h}Y_k)d\lambda\\
&+R_\tau D_{s}^{h}(W(t)-W(t_k))\\
&=D_{s}^{h}Y_k+(t-t_k)(B_\tau D_{s}^{h}Y_k+R_\tau G(Y_k)D_{s}^{h}Y_k).
\end{align*}
Therefore, since $|\tau B_\tau|_{\mathcal{L}(H)}\leq C$
$$|D_{s}^{h}\tilde{Y}(t)|_{\beta}\leq c|D_{s}^{h}Y_k|_{\beta},$$
and taking supremum over $h$ with $|h|\leq 1$ we get
\begin{equation}\label{mall}
|(-B)^\beta D_s\tilde{Y}(t)|_{\mathcal{L}(H)}\leq c|(-B)^\beta D_sY_k|_{\mathcal{L}(H)}.
\end{equation}
The last quantity is estimated thanks to Lemma \ref{lem5}:
$$|D_{s}^{h}Y_k|_\beta\leq C(1+L_G\tau)^{k-l_s}(1+\frac{1}{(1+\mu_0\tau)^{(1-\beta)(k-l_s)}t_{k-l_s}^{\beta}})|h|.$$
When $\tau\leq \tau_0$ and $(t_k-1)\vee 0\leq s\leq t_k\leq t< t_{k+1}$, we see that $(1+L_G\tau)^{k-l_s}$ is bounded by a constant.

We can then control the second term of $a_{k}^{1,3}$ with
\begin{align*}
&\tau\E\int_{t_{k}}^{t_{k+1}}\int_{(t_k-1)\vee 0}^{t_k}|R_\tau (-B)^{1/2+2\kappa}|_{\mathcal{L}(H)}|(-B)^{1-3\frac{\kappa}{2}}R_{\tau}^{k-l_s}|_{\mathcal{L}(H)}\text{Tr}((-B)^{-1/2-\frac{\kappa}{2}})\\
&\times|(-B)^{1/2-\kappa/2}D^2u(T-t,\tilde{Y}(t))(-B)^{1/2-\kappa/2}|_{\mathcal{L}(H)}|(-B)^\kappa D_s\tilde{Y}(t)|_{\mathcal{L}(H)}dsdt\\
&\leq C\tau^{1/2-2\kappa}\int_{t_{k}}^{t_{k+1}}\int_{(t_k-1)\vee 0}^{t_k}t_{k-l_s}^{-1+3\frac{\kappa}{2}}\frac{1}{(1+\mu_0\tau)^{(k-l_s)3\frac{\kappa}{2}}}(1+t_{k-l_s}^{-\kappa}\frac{1}{(1+\mu_0\tau)^{(k-l_s)(1-\kappa)}})ds\\
&\times(1+\frac{1}{(T-t)^\eta}+\frac{1}{(T-t)^{1-\kappa}})e^{-\tilde{\mu}(T-t)}(1+|y|^2)dt,
\end{align*}
using Proposition \ref{lem4} and Lemmas \ref{lem5} and \ref{lem6}.

On the one hand, we have

$$\int_{(t_k-1)\vee 0}^{t_k}t_{k-l_s}^{-1+3\frac{\kappa}{2}}\frac{1}{(1+\mu_0\tau)^{(k-l_s)3\frac{\kappa}{2}}}ds\leq\int_{0}^{t_k}\frac{1}{s^{1-3\frac{\kappa}{2}}}\frac{1}{(1+\mu_0\tau)^{3\frac{\kappa}{2} s/\tau}}ds\leq C<+\infty,$$
for $\tau\leq\tau_0$, thanks to \eqref{sumestim}.

On the other hand,
$$\sum_{k=1}^{m-1}\int_{t_k}^{t_{k+1}}(1+\frac{1}{(T-t)^\eta}+\frac{1}{(T-t)^{1-\kappa}})e^{-\tilde{\mu} (T-t)}dt\leq \int_{0}^{+\infty}(1+\frac{1}{t^\eta}+\frac{1}{t^{1-\kappa}})e^{-\tilde{\mu}t}dt<+\infty.$$

Therefore
\begin{equation}\label{e4}
\sum_{k=1}^{m-1}|a_{k}^{1,3}|\leq C(1+|y|^2)\tau^{1/2-2\kappa}.
\end{equation}

\end{enumerate}

\subsubsection{Estimate of $a_{k}^{2}$}

We decompose $a_{k}^{2}$ using the definition of $\tilde{Y}$ - see \eqref{tild}:
\begin{gather*}
a_{k}^{2,1}=\E\int_{t_k}^{t_{k+1}}(t-t_k)<BB_\tau Y_k,Du(T-t,\tilde{Y}(t))>dt\\
a_{k}^{2,2}=\E\int_{t_k}^{t_{k+1}}(t-t_k)<BR_\tau G(Y_k),Du(T-t,\tilde{Y}(t))>dt\\
a_{k}^{2,3}=\E\int_{t_k}^{t_{k+1}}<\int_{t_k}^{t}BR_\tau dW(s),Du(T-t,\tilde{Y}(t))>dt;
\end{gather*}
then $a_{k}^{2}=a_{k}^{2,1}+a_{k}^{2,2}+a_{k}^{2,3}$.

\begin{enumerate}
\item \textbf{Estimate of $a_{k}^{2,1}$}

Since $BB_\tau=R_\tau B^2$, $a_{k}^{2,1}$ is bounded by the same expression as $a_{k}^{1}$: by \eqref{e2}, \eqref{e3}, \eqref{e4} we have
\begin{equation}\label{e5}
\sum_{k=1}^{m-1}|a_{k}^{2,1}|\leq C(1+|y|^3)(1+T^{-(1/2-2\kappa)})\tau^{1/2-2\kappa}.
\end{equation}

\item \textbf{Estimate of $a_{k}^{2,2}$}

We have
\begin{align*}
|a_{k}^{2,2}|&\leq \tau \E\int_{t_{k}}^{t_{k+1}}|(-B)^{1/2+\kappa}R_{\tau}|_{\mathcal{L}(H)}|G(Y_k)||(-B)^{1/2-\kappa}Du(T-t,\tilde{Y}(t))|dt\\
&\leq \|G\|_\infty\tau^{1/2-\kappa}\int_{t_k}^{t_{k+1}}(1+\frac{1}{(T-t)^{1/2-\kappa}})e^{-\tilde{\mu}(T-t)}dt.
\end{align*}
We then have
\begin{equation}\label{e6}
\sum_{k=1}^{m-1}|a_{k}^{2,2}|\leq C\tau^{1/2-\kappa}.
\end{equation}

\item \textbf{Estimate of $a_{k}^{2,3}$}

We again use the integration by parts formula to rewrite $a_{k}^{2,3}$:
\begin{align*}
a_{k}^{2,3}&=\E\int_{t_k}^{t_{k+1}}<\int_{t_k}^{t}BR_\tau dW(s),Du(T-t,\tilde{Y}(t))>dt\\
&=\E\int_{t_{k}}^{t_{k+1}}\int_{t_k}^{t}\text{Tr}(R_\tau BD^2u(T-t,\tilde{Y}(t))D_s\tilde{Y}(t))dsdt.
\end{align*}
From \eqref{tild}, for $t_k\leq s\leq t\leq t_{k+1}$ we have $D_{s}^{h}\tilde{Y}(t)=R_\tau h$; as a consequence, we do not need to use the same trick as in the control of $a_{k}^{1,3}$.

Then we have
\begin{align*}
|a_{k}^{2,3}|&\leq\E\int_{t_k}^{t_{k+1}}(t-t_k)\text{Tr}(R_\tau B D^2u(T-t, \tilde{Y}(t))R_\tau)dt\\
&\leq c\tau\int_{t_{k}}^{t_{k+1}}|R_\tau(-B)^{1/2+\kappa/2}|_{\mathcal{L}(H)}\text{Tr}((-B)^{-1/2-\kappa/2})|(-B)^\kappa R_\tau|_{\mathcal{L}(H)}\\
&|(-B)^{1/2-\kappa/2}D^2u(T-t,\tilde{Y}(t))(-B)^{1/2-\kappa/2}|_{\mathcal{L}(H)}dt\\
&\leq c(1+|y|^2)\tau^{1/2-3\kappa/2}\int_{t_{k}}^{t_{k+1}}(1+\frac{1}{(T-t)^\eta}+\frac{1}{(T-t)^{1-\kappa}})e^{-\tilde{\mu}(T-t)}dt.
\end{align*}
Therefore
\begin{equation}\label{e7}
\sum_{k=1}^{m-1}|a_{k}^{2,3}|\leq C(1+|y|^2)\tau^{1/2-3\kappa/2}.
\end{equation}

\end{enumerate}

With the previous estimates on $a^{1}$ and $a^{2}$, we get
\begin{equation}\label{ea}
\sum_{k=1}^{m-1}|a_{k}|\leq C(1+|y|^3)(1+T^{-(1/2-2\kappa)})\tau^{1/2-2\kappa}.
\end{equation}

\subsection{Estimate of $b_k$}

We have
\begin{align*}
b_k&=\E\int_{t_k}^{t_{k+1}}<G(\tilde{Y}(t))-R_\tau G(Y_k),Du(T-t,\tilde{Y}(t))>dt\\
&=\E\int_{t_k}^{t_{k+1}}<(I-R_\tau)G(Y_k),Du(T-t,\tilde{Y}(t))>dt\\
&+\E\int_{t_k}^{t_{k+1}}<G(\tilde{Y}(t))-G(Y_k),Du(T-t,\tilde{Y}(t))>dt\\
&:=b_{k}^{1}+b_{k}^{2}.
\end{align*}

\subsubsection{Estimate of $b_{k}^{1}$}

This term is easy to treat: we have
\begin{align*}
|b_{k}^{1}|&\leq \E\int_{t_{k}}^{t_{k+1}}|(-B)^{-1/2+\kappa}(I-R_\tau)|_{\mathcal{L}(H)}|G(Y_k)||(-B)^{1/2-\kappa}Du(T-t,\tilde{Y}(t))|dt\\
&\leq C\tau^{1/2-\kappa}\int_{t_{k}}^{t_{k+1}}(1+\frac{1}{(T-t)^{1/2-\kappa}})e^{-\tilde{\mu}(T-t)}dt,
\end{align*}
where we have used Proposition \ref{lem3}, and the following inequality for $0\leq\beta\leq 1$:
\begin{equation}\label{normplus}
|(-B)^{-\beta}(I-R_\tau)|_{\mathcal{L}(H)}\leq C_\beta\tau^\beta.
\end{equation}
Then we see that
\begin{equation}\label{e8}
\sum_{k=1}^{m-1}|b_{k}^{1}|\leq C\tau^{1/2-\kappa}.
\end{equation}

\subsubsection{Estimate of $b_{k}^{2}$}

To estimate $|b_{k}^{2}|$, we write the scalar product in coordinates with respect to the orthonormal basis $(f_i)$, and then we expand the terms thanks to the It\^o formula.

If we note $G_i=<G,f_i>$ and $\partial_i=<D.,f_i>$, we have
$$<G(\tilde{Y}(t))-G(Y_k),Du(T-t,\tilde{Y}(t))>=\sum_{i}(G_i(\tilde{Y}(t))-G_i(Y_k))\partial_iu(T-t,\tilde{Y}(t)).$$

The above sum is finite, because we work with finite dimensional approximations.

It\^o formula gives for $t_k\leq t<t_{k+1}$
\begin{align*}
G_i(\tilde{Y}(t))-G_i(Y_k)&=\frac{1}{2}\int_{t_k}^{t}\text{Tr}(R_\tau R_{\tau}^{*}D^2G_i(\tilde{Y}(s)))ds\\
&+\int_{t_k}^{t}<B_\tau Y_k,DG_i(\tilde{Y}(s))>ds\\
&+\int_{t_k}^{t}<R_\tau G(Y_k),DG_i(\tilde{Y}(s))>ds\\
&+\int_{t_k}^{t}<DG_i(\tilde{Y}(s)),R_\tau dW(s)>.
\end{align*}

We naturally define $b_{k}^{2,j}$, for $j\in\left\{1,2,3,4\right\}$, and we now control each term.

\begin{enumerate}
\item \textbf{Estimate of $b_{k}^{2,1}$}

By definition, we have
$$b_{k}^{2,1}=\int_{t_k}^{t_{k+1}}\E\frac{1}{2}\int_{t_k}^{t}\sum_{i}\text{Tr}(R_\tau R_{\tau}^{*}D^2G_i(\tilde{Y}(s)))ds\partial_iu(T-t,\tilde{Y}(t))dt.$$
Using the orthonormal basis $(f_k)_k$ given by assumption \ref{hypB}, and recalling that the sums are finite, we can calculate:
\begin{align*}
\sum_{i}\text{Tr}(R_\tau R_{\tau}^{*}D^2G_i(\tilde{Y}(s)))&\partial_iu(T-t,\tilde{Y}(t))=\sum_{i}\text{Tr}(D^2G_i(\tilde{Y}(s))R_\tau R_{\tau}^{*})\partial_iu(T-t,\tilde{Y}(t))\\
&=\sum_{i}\sum_{j}<D^2G_i(\tilde{Y}(s))\frac{1}{(1+\mu_j\tau)^2}f_j,f_j>\partial_iu(T-t,\tilde{Y}(t))\\
&=\sum_{i}\sum_{j}\frac{1}{(1+\mu_j\tau)^2}D^2G_i(\tilde{Y}(s)).(f_j,f_j)\partial_iu(T-t,\tilde{Y}(t)).
\end{align*}
Using the Cauchy-Schwarz inequality (where $j$ is fixed), we get
\begin{multline*}
|\sum_{i}D^2G_i(\tilde{Y}(s)).(f_j,f_j)\partial_iu(T-t,\tilde{Y}(t))|\\
\leq \left(\sum_{i}\frac{|D^2G_i(\tilde{Y}(s)).(f_j,f_j)|^2}{\mu_{i}^{2\eta}}\right)^{1/2}\left(\sum_{i}\mu_{i}^{2\eta}|\partial_iu(T-t,\tilde{Y}(t))|^2\right)^{1/2}.
\end{multline*}
The second factor of this expression is $|(-B)^{\eta}Du(T-t,\tilde{Y}(t))|_{H}$; we control it thanks to Proposition \ref{lem3}. The first factor is controlled thanks to Assumption \ref{hypG}:
\begin{align*}
\left(\sum_{i}\frac{|D^2G_i(\tilde{Y}(s)).(f_j,f_j)|^2}{\mu_{i}^{2\eta}}\right)^{1/2}&=|(-B)^{-\eta}D^2G(\tilde{Y}(s)).(f_j,f_j)|\\
&\leq C|f_j|_H|f_j|_H\leq C,
\end{align*}
since $(f_j)_j$ is an orthonormal system.


Therefore
\begin{align*}
|\sum_{i}\text{Tr}(R_\tau R_{\tau}^{*}&D^2G_i(\tilde{Y}(s)))\partial_iu(T-t,\tilde{Y}(t))|\\
&\leq C(1+|y|^2)(1+\frac{1}{(T-t)^\eta})e^{-\tilde{\mu}(T-t)}\sum_{j=0}^{\infty}\frac{1}{(1+\mu_j\tau)^2}\\
&\leq C(1+|y|^2)(1+\frac{1}{(T-t)^\eta})e^{-\tilde{\mu}(T-t)} \tau^{-1/2-\kappa}\sum_{j=0}^{\infty}\frac{(\mu_j\tau)^{1/2+\kappa}}{(1+\mu_j\tau)^2}\frac{1}{\mu_{j}^{1/2+\kappa}}\\
&\leq C(1+|y|^2)(1+\frac{1}{(T-t)^\eta})e^{-\tilde{\mu}(T-t)}\tau^{-1/2-\kappa}.
\end{align*}
Then
$$|b_{k}^{2,1}|\leq C(1+|y|^2)\tau^{1/2-\kappa}\int_{t_{k}}^{t_{k+1}}(1+\frac{1}{(T-t)^\eta})e^{-\tilde{\mu}(T-t)}dt,$$
and
\begin{equation}\label{e9}
\sum_{k=1}^{m-1}|b_{k}^{2,1}|\leq C(1+|y|^2)\tau^{1/2-\kappa}.
\end{equation}

\item \textbf{Estimate of $b_{k}^{2,2}$}

Thanks to \eqref{exprYk} and \eqref{exprsto}, we have
\begin{align*}
b_{k}^{2,2}&=\E\int_{t_{k}}^{t_{k+1}}\int_{t_k}^{t}\sum_{i}<B_\tau R_{\tau}^{k}y+B_\tau\tau\sum_{l=0}^{k-1}R_{\tau}^{k-l}G(Y_l),DG_i(\tilde{Y}(s))>\partial_iu(T-t,\tilde{Y}(t))dsdt\\
&+\E\int_{t_{k}}^{t_{k+1}}\int_{t_k}^{t}\sum_{i}<B_\tau\int_{0}^{t_k}R_{\tau}^{k-l_r}dW(r),DG_i(\tilde{Y}(s))>\partial_iu(T-t,\tilde{Y}(t))dsdt\\
&:=b_{k}^{2,2,1}+b_{k}^{2,2,2}.
\end{align*}

{\em (i)} For the first term, recalling that $B_\tau=BR_\tau$ and that $G$ is bounded, we have
\begin{align*}
|b_{k}^{2,2,1}|&=|\E\int_{t_{k}}^{t_{k+1}}\int_{t_k}^{t}<DG(\tilde{Y}(s)).(B_\tau R_{\tau}^{k}y+B_\tau\tau\sum_{l=0}^{k-1}R_{\tau}^{k-l}G(Y_l)),Du(T-t,\tilde{Y}(t))>dsdt|\\
&\leq\E\int_{t_{k}}^{t_{k+1}}\int_{t_k}^{t}|(-B)^\kappa R_\tau|_{\mathcal{L}(H)}(|(-B)^{1-\kappa}R_{\tau}^{k}y|+\tau\sum_{l=0}^{k-1}|(-B)^{1-\kappa}R_{\tau}^{k-l}|_{\mathcal{L}(H)}|G(Y_l)|)\\
&\times|Du(T-t,\tilde{Y}(t))|dsdt\\
&\leq C\tau^{1-\kappa}\int_{t_k}^{t_{k+1}}(1+|y|^2)e^{-\tilde{\mu}(T-t)}dt(t_{k}^{-1+\kappa}|y|+\tau\sum_{l=0}^{k-1}t_{k-l}^{(1-\kappa)}\frac{1}{(1+\mu_0\tau)^{(k-l)\kappa}})\\
&\leq C\tau^{1-\kappa}(1+|y|^3)\frac{1}{t_{k}^{1-\kappa}}\int_{t_k}^{t_{k+1}}e^{-\tilde{\mu}(T-t)}dt,
\end{align*}
if $\tau\leq\tau_0$ - see \eqref{sumestim}.

Therefore
\begin{align*}
\sum_{k=1}^{m-1}|b_{k}^{2,2,1}|&\leq C\tau^{1-\kappa}(|y|+1)\sum_{k=1}^{m-1}\frac{1}{t_{k}^{1-\kappa}}\int_{t_k}^{t_{k+1}}e^{-\tilde{\mu}(T-t)}dt\\
&\leq C\tau^{1-\kappa}(1+|y|^3)\int_{0}^{T}\frac{1}{t^{1-\kappa}}e^{-\tilde{\mu}(T-t)}dt\\
&\leq C\tau^{1-\kappa}(1+|y|^3)\int_{0}^{T}\frac{1}{t^{1-\kappa}}\frac{C}{(T-t)^{1/2-\kappa}}dt\\
&\leq C\tau^{1-\kappa}(1+|y|^3)T^{-(1/2-2\kappa)}\int_{0}^{1}\frac{1}{s^{1-\kappa}}\frac{C}{(1-s)^{1/2-\kappa}}dt\\
&\leq C\tau^{1-\kappa}(1+|y|^3)T^{-(1/2-2\kappa)}.
\end{align*}

{\em (ii)} For the second term, we again use an integration by parts, after a decomposition of the time interval - as in the estimates for $a_{k}^{1,3}$.
First,
\begin{align*}
&b_{k}^{2,2,2}=\E\int_{t_{k}}^{t_{k+1}}\int_{t_k}^{t}\sum_{i}<B_\tau\int_{0}^{t_k}R_{\tau}^{k-l_r}dW(r),DG_i(\tilde{Y}(s))>\partial_iu(T-t,\tilde{Y}(t))dsdt\\
&=\E\int_{t_{k}}^{t_{k+1}}\int_{t_k}^{t}\sum_{i}<B_\tau\int_{0}^{(t_k-1)\vee 0}R_{\tau}^{k-l_r}dW(r),DG_i(\tilde{Y}(s))>\partial_iu(T-t,\tilde{Y}(t))dsdt\\
&+\E\int_{t_k}^{t_{k+1}}\int_{t_k}^{t}\sum_{i,j,m}<B_\tau\int_{(t_k-1)\vee 0}^{t_k}R_{\tau}^{k-l_r}f_m,f_j>d\beta_m(r)\partial_j G_i(\tilde{Y}(s))\partial_iu(T-t,\tilde{Y}(t))dsdt\\
&=:b_{k}^{2,2,2,1}+b_{k}^{2,2,2,2}.
\end{align*}


For $b_{k}^{2,2,2,1}$, we can work directly and see that
\begin{align*}
|b_{k}^{2,2,2,1}|&\leq |\E\int_{t_{k}}^{t_{k+1}}\int_{t_k}^{t}\sum_{i}<B_\tau\int_{0}^{(t_k-1)\vee 0}R_{\tau}^{k-l_r}dW(r),DG_i(\tilde{Y}(s))>\partial_iu(T-t,\tilde{Y}(t))dsdt|\\
&\leq \int_{t_k}^{t_{k+1}}\int_{t_k}^{t}\E|<DG(\tilde{Y}(s)).B_\tau\int_{0}^{(t_k-1)\vee 0}R_{\tau}^{k-l_r}dW(r),Du(T-t,\tilde{Y}(t))>|dsdt\\
&\leq \int_{t_k}^{t_{k+1}}\int_{t_k}^{t}(\E|B_\tau\int_{0}^{(t_k-1)\vee 0}R_{\tau}^{k-l_r}dW(r)|^{2})^{1/2}(\E|Du(T-t,\tilde{Y}(t))>|^2)^{1/2}dsdt\\
&\leq C\tau\int_{t_k}^{t_{k+1}}e^{-\tilde{\mu}(T-t)}(1+|y|^2),
\end{align*}
thanks to Lemmas \ref{lem2bis}, \ref{lem00} and to the following estimate for $\tau\leq \tau_0$
$$\E|B_\tau\int_{0}^{(t_k-1)\vee 0}R_{\tau}^{k-l_r}dW(r)|^{2}\leq\E|B^2R_\tau\int_{0}^{(t_k-1)\vee 0}R_{\tau}^{k-l_r}dW(r)|^{2}\leq C,$$
thanks to the estimate proved to control $a_{k}^{1,3}$.


For $b_{k}^{2,2,2,2}$, we can write thanks to a Malliavin integration by parts and with the chain rule
\begin{align*}
&b_{k}^{2,2,2,2}=\E\int_{t_k}^{t_{k+1}}\int_{t_k}^{t}\sum_{i,j,m}<B_\tau\int_{(t_k-1)\vee 0}^{t_k}R_{\tau}^{k-l_r}f_m,f_j>d\beta_m(r)\partial_j G_i(\tilde{Y}(s))\partial_iu(T-t,\tilde{Y}(t))dsdt\\
&=\E\int_{t_{k}}^{t_{k+1}}\int_{t_{k}}^{t}\int_{(t_k-1)\vee 0}^{t_k}\sum_{i,j,m,n}<B_\tau R_{\tau}^{k-l_r}f_m,f_j>\partial_{j,n}^{2}G_i(\tilde{Y}(s))<D_{r}^{m}\tilde{Y}(s),f_n>\partial_iu(T-t,\tilde{Y}(t))drdsdt\\
&+\E\int_{t_{k}}^{t_{k+1}}\int_{t_{k}}^{t}\int_{(t_k-1)\vee 0}^{t_k}\sum_{i,j,m,n}<B_\tau R_{\tau}^{k-l_r}f_m,f_j>\partial_jG_i(\tilde{Y}(s))\partial_{i,n}^{2}u(T-t,\tilde{Y}(t))<D_{r}^{m}\tilde{Y}(t),f_n>drdsdt\\
&=\E\int_{t_{k}}^{t_{k+1}}\int_{t_{k}}^{t}\int_{(t_k-1)\vee 0}^{t_k}\sum_{i,m}D^2G_i(\tilde{Y}(s))(B_\tau R_{\tau}^{k-l_r}f_m,D_{r}^{m}\tilde{Y}(s))\partial_iu(T-t,\tilde{Y}(t))drdsdt\\
&+\E\int_{t_{k}}^{t_{k+1}}\int_{t_{k}}^{t}\int_{(t_k-1)\vee 0}^{t_k}\sum_{i,m}<\mathcal{B}_i(s,t)B_\tau R_{\tau}^{k-l_r}f_m,D_{r}^{m}\tilde{Y}(t)>drdsdt\\
&=\E\int_{t_{k}}^{t_{k+1}}\int_{t_{k}}^{t}\int_{(t_k-1)\vee 0}^{t_k}\sum_{i}\text{Tr}\left((D_{r}\tilde{Y}(s))^{*}D^2G_i(\tilde{Y}(s))B_\tau R_{\tau}^{k-l_r}\right)\partial_iu(T-t,\tilde{Y}(t))drdsdt\\
&+\E\int_{t_{k}}^{t_{k+1}}\int_{t_{k}}^{t}\int_{(t_k-1)\vee 0}^{t_k}\sum_{i}\text{Tr}\left((D_r\tilde{Y}(t)^{*}\mathcal{B}_i(s,t)B_\tau R_{\tau}^{k-l_r}\right)drdsdt,
\end{align*}
where we define a linear operator on $H$ by
\begin{align*}
<\mathcal{B}_i(s,t)h,k>&=<DG_i(\tilde{Y}(s)),h>\sum_{n=0}^{+\infty}\partial_{i,n}^{2}u(T-t,\tilde{Y}(t))<k,f_n>\\
&=<DG_i(\tilde{Y}(s)),h><D^2u(T-t,\tilde{Y}(t)).f_i,k>.
\end{align*}

We have $\sum_{i}<\mathcal{B}_i(s,t)h,k>=D^2u(T-t,\tilde{Y}(t)).(DG(\tilde{Y}(s)).h,k)$,
and $$|\sum_{i}\mathcal{B}_i(s,t)|_{\mathcal{L}(H)}\leq|DG(\tilde{Y}(s))|_{\mathcal{L}(H)}|D^2u(T-t,\tilde{Y}(t))|_{\mathcal{L}(H)};$$
so we can write, for $(t_k-1)\vee 0\leq r\leq t_k$
\begin{align*}
&|\sum_i\text{Tr}\left((D_r\tilde{Y}(t)^{*}\mathcal{B}_i(s,t)B_\tau R_{\tau}^{k-l_r}\right)|\\
&\leq|D_r\tilde{Y}(t)|_{\mathcal{L}(H)}|\sum_{i}\mathcal{B}_i(s,t)|_{\mathcal{L}(H)}|(-B)^{1-3\kappa/2}R_{\tau}^{k-l_r}|_{\mathcal{L}(H)}|R_\tau(-B)^{1/2+2\kappa}|_{\mathcal{L}(H)}\text{Tr}((-B)^{-1/2-\kappa/2})\\
&\leq C\tau^{-1/2-2\kappa}t_{k-l_r}^{-1+3\kappa/2}\frac{1}{(1+\mu_0\tau)^{(k-l_r)3\kappa/2}}e^{-\tilde{\mu}(T-t)},
\end{align*}
using Proposition \ref{lem4}, Lemma \ref{lem5} - since $(1+L_G\tau)^{k-l_r}\leq C$ - Lemma \ref{lem6} and estimate \eqref{mall}.

The other term is a little more complicated, because we are not able to control $D^2G(\tilde{Y}(s))$ in $H$. We proceed as in the estimate of $b_{k}^{2,1}$, and we directly calculate the trace.
\begin{align*}
&|\sum_{i}\text{Tr}\left((D_{r}\tilde{Y}(s))^{*}D^2G_i(\tilde{Y}(s))B_\tau R_{\tau}^{k-l_r}\right)\partial_iu(T-t,\tilde{Y}(t))|\\
&\leq |D_r\tilde{Y}(s)|_{\mathcal{L}(H)}|\sum_{i}\text{Tr}\left(D^2G_i(\tilde{Y}(s))B_\tau R_{\tau}^{k-l_r}\right)\partial_iu(T-t,\tilde{Y}(t))|\\
&\leq |D_r\tilde{Y}(s)|_{\mathcal{L}(H)}\sum_{i,j}\frac{|D^2G_i(\tilde{Y}(s)).(f_j,f_j)|}{\mu_{i}^{\eta}}\frac{\mu_j}{(1+\mu_j\tau)^{1+k-l_r}}\mu_{i}^{\eta}|\partial_iu(T-t,\tilde{Y}(t))|\\
&\leq |D_r\tilde{Y}(s)|_{\mathcal{L}(H)}|(-B)^\eta Du(T-t,\tilde{Y}(t))|_{H}\sum_{j}|(-B)^{-\eta}D^2G(\tilde{Y}(s)).(f_j,f_j)|\frac{\mu_j}{(1+\mu_j\tau)^{1+k-l_r}},
\end{align*}
thanks to the Cauchy-Schwarz inequality.

By using the same analysis as in the estimation of $b_{k}^{2,1}$, we see that the above expression is bounded by
$$C|D_r\tilde{Y}(s)|_{\mathcal{L}(H)}|(-B)^\eta Du(T-t,\tilde{Y}(t))|_{H}\sum_{j}\frac{\mu_j}{(1+\mu_j\tau)^{1+k-l_r}};$$
but the last sum is equal to $\text{Tr}(B_\tau R_{\tau}^{k-l_r})$, so that we see that indeed the two expressions in $b_{k}^{2,2,2}$ are bounded by the same expression.

Therefore
\begin{align*}
&|b_{k}^{2,2,2,2}|\\
&\leq \E\int_{t_{k}}^{t_{k+1}}\int_{t_{k}}^{t}\int_{(t_k-1)\vee 0}^{t_k}C\tau^{-1/2-2\kappa}t_{k-l_r}^{-1+3\kappa/2}\frac{e^{-\tilde{\mu}(T-t)}}{(1+\mu_0\tau)^{(k-l_r)3\kappa/2}}(1+\frac{1}{(T-t)^\eta})(1+|y|^2)drdsdt\\
&\leq C(1+|y|^2)\tau^{1/2-2\kappa}\int_{t_{k}}^{t_{k+1}}(1+\frac{1}{(T-t)^\eta})e^{-\tilde{\mu}(T-t)}dt\int_{0}^{t_k}t_{k-l_r}^{-1+3\kappa/2}\frac{1}{(1+\mu_0\tau)^{(k-l_r)3\kappa/2}}dr\\
&\leq C(1+|y|^2)\tau^{1/2-2\kappa}\int_{t_{k}}^{t_{k+1}}(1+\frac{1}{(T-t)^\eta})e^{-\tilde{\mu}(T-t)}dt,
\end{align*}
as already proved - see \eqref{sumestim}.

Now gathering estimates for $b_{k}^{2,2,2,1}$ and $b_{k}^{2,2,2,2}$, we obtain
\begin{equation}\label{e10}
\sum_{k=1}^{m-1}|b_{k}^{2,2,2}|\leq C(1+|y|^2)\tau^{1/2-2\kappa}.
\end{equation}

\item \textbf{Estimate of $b_{k}^{2,3}$}
We have
\begin{align*}
b_{k}^{2,3}&=\E\int_{t_{k}}^{t_{k+1}}\int_{t_k}^{t}\sum_{i}<R_\tau G(Y_k),DG_i(\tilde{Y}(s))>\partial_iu(T-t,\tilde{Y}(t))dsdt\\
&=\E\int_{t_{k}}^{t_{k+1}}\int_{t_k}^{t}<Du(T-t,\tilde{Y}(t)),DG(\tilde{Y}(s)).(R_\tau G(Y_k))>dsdt.
\end{align*}
Using that $G$ and $DG$ are bounded, we easily see that
$$|b_{k}^{2,3}|\leq C(1+|y|^2)\tau\int_{t_k}^{t_{k+1}}e^{-\tilde{\mu}(T-t)}dt,$$
and that
\begin{equation}\label{e11}
\sum_{k=1}^{m-1}|b_{k}^{2,3}|\leq C(1+|y|^2)\tau.
\end{equation}

\item \textbf{Estimate of $b_{k}^{2,4}$}

We use the integration by parts formula of Proposition \ref{lemintbyparts} to get
\begin{align*}
b_{k}^{2,4}&=\int_{t_{k}}^{t_{k+1}}\int_{t_k}^{t}\sum_{i}<DG_i(\tilde{Y}(s)),R_\tau dW(s)>\partial_iu(T-t,\tilde{Y}(t))dt\\
&=\E\int_{t_k}^{t_{k+1}}\int_{t_k}^{t}\text{Tr}\left((D_s\tilde{Y}(t))^{*}D^2u(T-t,\tilde{Y}(t))DG(\tilde{Y}(s))R_\tau\right)dsdt\\
&=\E\int_{t_k}^{t_{k+1}}\int_{t_k}^{t}\text{Tr}\left(R_\tau D^2u(T-t,\tilde{Y}(t))DG(\tilde{Y}(s))R_\tau\right)dsdt,
\end{align*}
using the identity $D_{s}^{h}\tilde{Y}(t)=R_\tau h$ when $t_k\leq s\leq t\leq t_{k+1}$, as in the estimate of $a_{k}^{2,3}$.

Now
\begin{align*}
|b_{k}^{2,4}|&\leq\E\int_{t_k}^{t_{k+1}}\int_{t_k}^{t}|(R_\tau(-B)^{1/2+\kappa}|_{\mathcal{L}(H)}|DG(\tilde{Y}(s))|_{\mathcal{L}(H)}|R_\tau|_{\mathcal{L}(H)}\\
&\times |D^2u(T-t,\tilde{Y}(t))|_{\mathcal{L}(H)}\text{Tr}((-B)^{-1/2-\kappa})dsdt\\
&\leq C(1+|y|^2)\tau^{1/2-\kappa}\int_{t_k}^{t_{k+1}}(1+\frac{1}{(T-t)^{\eta}})e^{-\tilde{\mu}(T-t)}dt,
\end{align*}
and
\begin{equation}\label{e12}
\sum_{k=1}^{m-1}|b_{k}^{2,4}|\leq  C(1+|y|^2)\tau^{1/2-\kappa}.
\end{equation}

\end{enumerate}

\subsubsection{Estimate of $b_k$: conclusion}
With \eqref{e8}, \eqref{e9}, \eqref{e10}, \eqref{e11} and \eqref{e12}, we get
\begin{equation}\label{eb}
\sum_{k=1}^{m-1}|b_k|\leq C\tau^{1/2-2\kappa}.
\end{equation}

\subsection{Estimate of $c_k$}

We have, using the symmetry of $R_\tau$,
$$\frac{1}{2}I-\frac{1}{2}R_\tau R_{\tau}^{*}=R_\tau(I-R_\tau)^*+\frac{1}{2}(I-R_\tau)(I-R_\tau)^*,$$
and
\begin{align*}
c_k&=\frac{1}{2}\E\int_{t_k}^{t_{k+1}}\text{Tr}((I-R_\tau R_{\tau}^{*})D^2u(T-t,\tilde{Y}(t)))dt\\
&=\frac{1}{2}\E\int_{t_k}^{t_{k+1}}\text{Tr}((I-R_\tau)(I-R_\tau)^{*}D^2u(T-t,\tilde{Y}(t)))dt\\
&+\E\int_{t_k}^{t_{k+1}}\text{Tr}(R_\tau(I-R_\tau)^{*}D^2u(T-t,\tilde{Y}(t)))dt\\
&:=c_{k}^{1}+c_{k}^{2}.
\end{align*}

\subsubsection{Estimate of $c_{k}^{1}$}

We have, using inequality \eqref{normplus}
\begin{align*}
|c_{k}^{1}|&\leq \frac{1}{2}\E\int_{t_k}^{t_{k+1}}\text{Tr}((-B)^{-1/2+\kappa}(I-R_\tau)^2(-B)^{-1/2+\kappa})\\
&\hspace{50 pt}\times|(-B)^{1/2-\kappa}D^2u(T-t,\tilde{Y}(t))(-B)^{1/2-\kappa}|_{\mathcal{L}(H)}dt\\
&\leq C(1+|y|^2)\int_{t_{k}}^{t_{k+1}}|(-B)^{-1/2+3\kappa}(I-R_\tau)|_{\mathcal{L}(H)}|I-R_\tau|_{\mathcal{L}(H)}\text{Tr}((-B)^{-1/2-\kappa})\\
&\hspace{50 pt}\times(1+\frac{1}{(T-t)^{\eta}}+\frac{1}{(T-t)^{1-\kappa}})e^{-\tilde{\mu}(T-t)}dt\\
&\leq C(1+|y|^2)\tau^{1/2-3\kappa}\int_{t_{k}}^{t_{k+1}}(1+\frac{1}{(T-t)^{\eta}}+\frac{1}{(T-t)^{1-\kappa}})e^{-\tilde{\mu}(T-t)}dt.
\end{align*}

Then
\begin{equation}\label{e13}
\sum_{k=1}^{m-1}|c_{k}^{1}|\leq C(1+|y|^2)\tau^{1/2-3\kappa}.
\end{equation}

\subsubsection{Estimate of $c_{k}^{2}$}

We have, using inequality \eqref{normplus}
\begin{align*}
|c_{k}^{2}|&\leq \E\int_{t_{k}}^{t_{k+1}}\text{Tr}((-B)^{-1/2+\kappa}R_\tau(I-R_\tau)(-B)^{-1/2+\kappa})\\
&\hspace{50 pt}\times|(-B)^{1/2-\kappa}D^2u(T-t,\tilde{Y}(t))(-B)^{1/2-\kappa}|_{\mathcal{L}(H)}dt\\
&\leq C(1+|y|^2)\int_{t_{k}}^{t_{k+1}}|(-B)^{-1/2+\kappa}(I-R_\tau)(-B)^{2\kappa}|_{\mathcal{L}(H)}\text{Tr}((-B)^{-1/2-\kappa})\\
&\hspace{50 pt}\times(1+\frac{1}{(T-t)^{\eta}}+\frac{1}{(T-t)^{1-\kappa}})e^{-\tilde{\mu}(T-t)}dt\\
&\leq C(1+|y|^2)\tau^{1/2-3\kappa}\int_{t_{k}}^{t_{k+1}}(1+\frac{1}{(T-t)^{\eta}}+\frac{1}{(T-t)^{1-\kappa}})e^{-\tilde{\mu}(T-t)}dt.
\end{align*}

Then
\begin{equation}\label{e14}
\sum_{k=1}^{m-1}|c_{k}^{2}|\leq C(1+|y|^2)\tau^{1/2-3\kappa}.
\end{equation}

\subsubsection{Estimate of $c_k$: conclusion}
With \eqref{e13} and \eqref{e14}, we get
\begin{equation}\label{ec}
\sum_{k=1}^{m-1}|c_k|\leq C(1+|y|^2)\tau^{1/2-3\kappa}.
\end{equation}

\subsection{Conclusion}

We put together estimates \eqref{ea}, \eqref{eb}, \eqref{ec} and \eqref{e1}; then passing to the limit with respect to dimension, we get the result.


\end{document}